\documentclass[11pt]{amsart}

\usepackage{amsmath, amsthm, amssymb, fullpage}

\usepackage{paralist}
\usepackage{enumerate}

\usepackage{hyperref} 

\newcommand{\rd}{\mathrm{d}}
\newcommand{\rM}{\mathrm{M}}
\newcommand{\s}{\mathrm{s}}
\renewcommand{\ss}{\mathrm{ss}}
\newcommand{\cR}{\mathcal{R}}
\newcommand{\cP}{\mathcal{P}}
\newcommand{\bA}{\mathbb{A}}
\newcommand{\bP}{\mathbb{P}}  
\newcommand{\bR}{\mathbb{R}}
\newcommand{\bC}{\mathbb{C}} 
\newcommand{\sP}{\mathsf{P}}
\newcommand{\sH}{\mathsf{H}}
\DeclareMathOperator{\spe}{sp}
\DeclareMathOperator{\ord}{ord} 
\DeclareMathOperator{\PGL}{PGL}
\DeclareMathOperator{\dep}{depth} 
\DeclareMathOperator{\hypRes}{hypRes}
\DeclareMathOperator{\ordRes}{ordRes}
\DeclareMathOperator{\MHRL}{MHRL}
\DeclareMathOperator{\supp}{Supp}
\DeclareMathOperator{\Res}{Res}
\DeclareMathOperator{\Spec}{Spec}
\DeclareMathOperator{\Proj}{Proj}
\DeclareMathOperator{\Rat}{Rat}
\DeclareMathOperator{\IRD}{IRD}
\DeclareMathOperator{\SL}{SL}
\DeclareMathOperator{\Sym}{Sym}
\newcommand{\Vu}{\boldsymbol{u}}
\newcommand{\Vv}{\boldsymbol{v}}
\newcommand{\Vw}{\boldsymbol{w}}
\newcommand{\widevec}[1]{\overrightarrow{#1}}

\theoremstyle{plain}
\newtheorem{theorem}{Theorem}[section]
\newtheorem{lemma}[theorem]{Lemma}
\newtheorem{proposition}[theorem]{Proposition}

\newtheorem{mainth}{Theorem}
\newtheorem{maincoro}[mainth]{Corollary}

\theoremstyle{definition}

\newtheorem*{acknowledgement}{Acknowledgement}

\theoremstyle{remark}

\theoremstyle{plain}
\newtheorem{apptheorem}{Theorem}[section]
\newtheorem{applemma}[apptheorem]{Lemma}

\theoremstyle{remark}

\theoremstyle{definition}

\numberwithin{equation}{section}

\makeatletter
\@namedef{subjclassname@2020}{\textup{2020} Mathematics Subject Classification}
\makeatother

\begin{document}

\title[Intrinsic reduction data and indeterminacy loci]
{Intrinsic reduction data and the loci related to indeterminacy of iterate maps for rational functions}
\author{Y\^usuke Okuyama}

\address{Division of Mathematics, Kyoto Institute of Technology, Sakyo-ku, Kyoto 606-8585 JAPAN}
\email{okuyama@kit.ac.jp}

\keywords{dynamical moduli, intrinsic reduction data, 
iterate maps, indeterminacy locus, dimension, Berkovich projective line}
\subjclass[2020]{Primary 37P45; Secondary 14L24, 37P50}

\date{\today}

\begin{abstract}
Let $k$ be an algebraically closed field of not necessarily characteristic $0$, and $d,n$ be integers $\ge 2$. We first formulate the dynamical moduli spaces of rational functions defined over $k$ in terms of the intrinsic reduction data of non-archimedean rational functions.   
As an application, we show that the locus $C_{d,n}$ of semistable intrinsic reduction data in degree $d$ whose $n$th iteration is defined and unstable is of codimension $d$ in the intrinsic reduction data space of degree $d$, settling a conjecture of Kiwi and Nie, and 
show that the image of $C_{d,n}$ in the compactified dynamical moduli $\overline{\rM}_{d,k}$ is of dimension $\le d-2$ and the equality holds precisely either if $d$ is even or if $d\ge 5$ is odd and $(4/d)^n<1/(d-3)$ and determine the dimensions even for all the five exceptional pairs $(d,n)$, 
in addition describing the images of $C_{3,n}$ themselves.
We also determine the indeterminacy of iterate maps  
either on the intrinsic reduction data spaces or on the compactified dynamical moduli spaces.
\end{abstract}
\maketitle

\section{Introduction}\label{sec:intro}

\subsection{Intrinsic formulation of dynamical moduli}\label{sec:form}

For any ring $R$ and any integer $e\ge 0$, set $V_{e,R}=R^2\otimes_R\Sym^e\bigl((R^2)^\vee\bigr)$. When $R$ is any field $K$,
for each $(F_0,F_1)\in\Sym^e((K^2)^\vee)^2\setminus\{(0,0)\}$, 
the projective class of $(1,0)\otimes F_0+(0,1)\otimes F_1$ in $\bP(V_{e,K})$ is denoted by $[F_0:F_1]$.
For every $e\ge 0$, letting $\Rat_{e,K}$ denote the set of all rational functions of degree $e$ over $K$, we have the bijection
\begin{gather*}
 \Rat_{e,K}\to \bigl(D_+(\Res)\bigr)(K),\quad 
P/Q\mapsto\bigl[X^eQ(Y/X):X^eP(Y/X)\bigr],
\end{gather*}
where $P,Q\in K[z]$, $\gcd(P,Q)\in K^\times$, and 
$\max\{\deg P,\deg Q\}=e$;
$\Rat_{e,K}$ also denotes the open subvariety $D_+(\Res)$ of $\bP(V_{e,K})$.
For any valuation ring $R$, let $\eta_R$ denote the generic point of
$\Spec R$, and let $\spe_R$ be the reduction map that sends each section
over $\eta_R$ of a proper $R$-scheme to the value at the closed point in
$\Spec R$ of the unique extension to $\Spec R$ of that section.

Now let $K$ be an algebraically closed field that is complete with respect to
a non-trivial and non-archimedean absolute value, and
$k=K^\circ/K^{\circ\circ}$ the residue field of $K$, where
$K^\circ$ and $K^{\circ\circ}$ denote
the valuation ring of $K$ and the maximal ideal of $K^\circ$, respectively. For every $d\ge 0$, by properness of $\bP(V_{d,K^\circ})$, 
we have the reduction map $\spe_{K^\circ}:\bP(V_{d,K})(K)\to\bP(V_{d,k})(k)$. The restriction to $\Rat_{d,K}$ of the reduction map 
is surjective onto $\bP(V_{d,k})$, and we denote by $\widehat{\phi}=\spe_{K^\circ}(\phi)$ the image of each $\phi\in\Rat_{d,K}$ and call it the coefficient reduction of $\phi$. The action on $\bP^1$ of each $\phi\in\Rat_{d,K}$ 
canonically extends to that on the Berkovich projective line
$\sP^1=\sP^1_K$. When $d>0$,
this canonical action is still $d$ to $1$ taking into account the local degree $\deg_\xi(\phi)$ of $\phi$ at each $\xi\in\sP^1$ and,
for every $\xi\in\sP^1$, 
induces the tangent map $\phi_{*,\xi}:T_\xi\sP^1\to T_{\phi(\xi)}\sP^1$ 
of $\phi$ at $\xi$ between the tangent (or direction) spaces of $\sP^1$
at $\xi,\phi(\xi)$, and moreover, the intrinsic reduction
$\tilde{\phi}_\xi:T_\xi\sP^1\to T_\xi\sP^1$ of $\phi$ at $\xi\in\sP^1$
and the intrinsic depth of $\phi$ at each direction $\Vv\in T_\xi\sP^1$ are defined by
\begin{gather}
 \tilde{\phi}_\xi
 \begin{cases}
 :=\phi_{*,\xi}&\text{if }\phi(\xi)=\xi,\\
 :\equiv\widevec{\xi(\phi(\xi))}&\text{otherwise}
 \end{cases}
 \quad\text{and}\quad
 \dep_{\Vv}\tilde{\phi}_\xi:=
 (\phi^*\delta_\xi)(U(\Vv))\in\{0,\ldots,d\},\label{eq:intrinsic}
\end{gather}
so that for every $\xi\in\sP^1$,
$\delta_\xi$ is the Dirac measure on $\sP^1$ at $\xi$ and
$\phi^*\delta_\xi:=\sum_{\xi'\in\phi^{-1}(\xi)}\deg_{\xi'}(\phi)\delta_{\xi'}$. The ordered pair
$(\tilde{\phi}_\xi,\dep_{\,\cdot}\tilde{\phi}_{\xi})$
is called the intrinsic reduction data of $\phi\in\Rat_{d,K}$ at each $\xi\in\sP^1$ and
we denote by $\IRD_{d,k}(\xi)$ the totality of all such intrinsic reduction data $(\tilde{\phi}_\xi,\dep_{\,\cdot}\tilde{\phi}_{\xi})$,
$\phi\in\Rat_{d,K}$.
When $\xi$ is type II and is fixed by $\phi$,  
the intrinsic reduction $\tilde{\phi}_\xi$ and the intrinsic depth function $\dep_{\,\cdot}\tilde{\phi}_{\xi}$ on $T_{\xi}\sP^1$
are regarded as an element of $\Rat_{\deg_\xi(\phi),k}$ 
and an effective divisor on $\bP^1_k$ of total $d-\deg_\xi(\phi)$, 
with respect to a natural $k$-projective structure on $T_\xi\sP^1$.

For every $d>0$, denoting by $\xi_g$ 
the Gauss (or canonical type II) point in $\sP^1$,
there is a $\PGL(2,K^\circ)$-equivariant bijection
\begin{gather}
 \IRD_{d,k}(\xi_g)\to\bP(V_{d,k}),
 \quad(\tilde{\phi}_{\xi_g},\dep_{\,\cdot}\tilde{\phi}_{\xi_g})\mapsto\hat{\phi},\label{eq:ird}
\end{gather}
where writing $\hat{\phi}=[HG_0:HG_1]$, 
$e\in\{0,\ldots,d\}$,
$H\in\Sym^{d-e}((k^2)^\vee)\setminus\{0\}$, and 
$G_0,G_1\in\Sym^e((k^2)^\vee)$ such that $\gcd(G_0,G_1)\in k^\times$,
$\tilde{\phi}_{\xi_g}=[G_0:G_1]$ and $\dep_{\,\cdot}\tilde{\phi}_{\xi_g}=\sum_{\hat{a}\in\bP^1(k)}(\ord_{\hat{a}}[H=0])1_{\{\widevec{\xi_ga}\}}$.
This bijection induces a canonical algebraic structure of $\IRD_{d,k}(\xi_g)$ in one direction and gives a non-archimedean interpretation of $\bP(V_{d,k})$ in the other direction, and identifies $\Rat_{d,k}\subset\bP(V_{d,k})$ with 
the open subvariety $\{(\tilde{\phi}_{\xi_g},0):\tilde{\phi}_{\xi_g}\in\Rat_{d,k}\}$ in $\IRD_{d,k}(\xi_g)$.

For every $d>1$, let us formulate the dynamical moduli space $\rM_{d,k}$
and the compactified dynamical moduli space $\overline{\rM}_{d,k}$
of degree $d$ from GIT in terms of intrinsic reduction data at $\xi_g$.
The conjugation action on $\bP(V_{d,k})$
of $\PGL_{2,k}=\Rat_{1,k}=D_+(\Res)\subset\bP(V_{1,k})$
is identified with that on $\IRD_{d,k}(\xi_g)$
of $\PGL(2,K^\circ)$ or 
strictly speaking of
$\PGL(2,K^\circ)/\bigl\{M\in\PGL(2,K^\circ):\tilde{M}_{\xi_g}=\mathrm{id}\bigr\}$. The GIT semistable and stable loci
$\IRD_{d,k}^{\ss}(\xi_g)$ and $\IRD_{d,k}^{\s}(\xi_g)$
under this conjugation action of $\PGL(2,K^\circ)$
are non-empty open subsets of $\IRD_{d,k}(\xi_g)$,
and $\IRD_{d,k}^{\s}(\xi_g)$ is a saturated open subset of
$\IRD_{d,k}^{\ss}(\xi_g)$ under
the categorical quotient morphism
\begin{gather*}
 q_{d,k}:\IRD_{d,k}^{\ss}(\xi_g)\to\overline{\rM}_{d,k}
 :=\IRD_{d,k}^{\ss}(\xi_g)//\PGL(2,K^\circ),
\end{gather*}
which restricts to the geometric quotient morphisms
\begin{gather*}
 \IRD_{d,k}^{\s}(\xi_g)\to\IRD_{d,k}^{\s}(\xi_g)/\SL_{2,k}
 \quad\text{and}\quad  \Rat_{d,k}\to\rM_{d,k}:=\Rat_{d,k}/\SL_{2,k}
\end{gather*}
(for some details on this GIT-semistability, see Proof of Theorem \ref{th:coeff}).
Those GIT-(semi)stability notions
are also intrinsically characterized or formulated as follows
(\cite[\S1.1 and Theorem 2.1]{OkuyamaDepth}); for every $\phi\in\Rat_{d,K}$, 
\begin{gather}
 \begin{cases}
 (\tilde{\phi}_{\xi_g},\dep_{\,\cdot}\tilde{\phi}_{\xi_g})\in\IRD_{d,k}^{\ss}(\xi_g)
 \Leftrightarrow \xi_g\in\MHRL(\phi),\\
 (\tilde{\phi}_{\xi_g},\dep_{\,\cdot}\tilde{\phi}_{\xi_g})\in\IRD_{d,k}^{\s}(\xi_g)
 \Leftrightarrow \MHRL(\phi)=\{\xi_g\},
 \end{cases}\label{eq:git}
\end{gather}
where $\MHRL(\phi)$ is the minimum locus of the hyperbolic resultant function
\begin{gather*}
 \hypRes_\phi(\xi):=
 \frac{\rho(\xi,\xi_g)}{2}
 +\frac{\rho(\xi,\phi(\xi)\wedge_{\xi_g}\xi)
       -\int_{\sP^1}\rho(\xi_g,\xi\wedge_{\xi_g}\cdot)
                         \,\rd(\phi^*\delta_{\xi_g})}{d-1}
\end{gather*}
on $\sP^1$ setting $\equiv+\infty$ on $\bP^1$. For more details on
the above intrinsic notions, see Section \ref{sec:background}.

For every type II point $\xi\in\MHRL(\phi)$ and every
$M\in\PGL(2,K)$ such that $M(\xi)=\xi_g$, the image
$q_{d,k}(\widetilde{(M\circ\phi\circ M^{-1})}_{\xi_g},\dep_{\,\cdot}\widetilde{(M\circ\phi\circ M^{-1})}_{\xi_g})
\in\overline{\rM}_{d,k}$ of the intrinsic reduction data
is independent of $\xi$ and $M$; see, e.g.,
\cite[Fact 5.1(i)]{OkuyamaQuadratic}.

\subsection{Dimensions of loci related to indeterminacy of iterate maps}\label{sec:mainresults}

For every $d>1$ and every $n\ge 2$, the $n$th iterate map
\begin{gather*}
 \Phi^{\mathrm{IRD}}_{d,n,k}:\IRD_{d,k}(\xi_g)\dashrightarrow\IRD_{d^n,k}(\xi_g)
\end{gather*}
on $\IRD_{d,k}(\xi_g)$
is a rational extension of the $n$th iterate morphism $\Rat_{d,k}\to\Rat_{d^n,k}$ and the indeterminacy locus $I(\Phi^{\mathrm{IRD}}_{d,n,k})$
of $\Phi^{\mathrm{IRD}}_{d,n,k}$
is independent of $n$ and equals the locus
\begin{gather*}
 I(d):=\Bigl\{(\tilde{\phi}_{\xi_g},\dep_{\,\cdot}\tilde{\phi}_{\xi_g})\in\IRD_{d,k}(\xi_g):\phi(\xi_g)\neq\xi_g
\text{ and }\dep_{\widevec{\xi_g\phi(\xi_g)}}\tilde{\phi}_{\xi_g}>0
\Bigr\}
\end{gather*}
in $\IRD_{d,k}(\xi_g)$; we give a proof of this fact in Theorem \ref{th:coeff} (cf.\ DeMarco \cite[Lemma 2.2]{DeMarco}).
For every $n\ge 2$, the locus
\begin{gather}
 C_{d,n}:=
 \bigl\{(\tilde{\phi}_{\xi_g},\dep_{\,\cdot}\tilde{\phi}_{\xi_g})\in\IRD_{d,k}(\xi_g)\setminus I(d):\xi_g\in\MHRL(\phi)\setminus\MHRL(\phi^n)\bigr\}\label{eq:unstable}
\end{gather}
is a closed subset in
$\bigl\{(\tilde{\phi}_{\xi_g},\dep_{\,\cdot}\tilde{\phi}_{\xi_g})
\in\IRD_{d,k}(\xi_g)\setminus I(d):\xi_g\in\MHRL(\phi)\bigr\}$, so that
the image $q_{d,k}(C_{d,n})$ is constructible by Chevalley's theorem
(see \cite[II, Exercise 3.19]{Hartshorne}). 
Here and below, the dimension of a constructible set means that
of the Zariski closure of this set.
When $k=\bC$, the locus $C_{d,n}$ coincides with Kiwi--Nie's $U_n=U_{d,n}$ in \cite{KN23},
where $d>1$ is implicit. 

They expected in \cite[Remark 3.5]{KN23} that the equalities
$\operatorname{codim}_{\bP^{2d+1}_k}U_n=d$ and
$\dim q_{d,k}(U_n)=d-2$ would hold for every $n\ge 2$. 
It turns out that our intrinsic formulation above is very effective
in computing $C_{d,n}$ (so $U_{d,n}$) and $I(d)$ in $\IRD_{d,k}(\xi_g)$ and their images 
in $\overline{\rM}_{d,k}$ under $q_{d,k}$.
Our first principal result is the following, which answers the above
question of Kiwi--Nie in a more precise form than they expected.

\begin{mainth}\label{th:main} 
Let $d\ge 2$ and $n\ge 2$.
\begin{inparaenum}[(i)]
\item\label{item:coeff} We have the equalities
\begin{gather}
 \dim C_{d,n}=d+1
 \quad\text{and}\quad
 \operatorname{codim}_{\IRD_{d,k}(\xi_g)}C_{d,n}=d.\label{eq:main}
\end{gather}
\item\label{item:moduli} In general, we have the inequality
$\dim q_{d,k}(C_{d,n})\le d-2$.
When $d$ is even,
\begin{gather}
 \dim q_{d,k}(C_{d,n})=d-2.\label{eq:evenimage}
\end{gather}
When $d$ is odd and $\ge 5$, \eqref{eq:evenimage} holds if and only if
$(4/d)^n<1/(d-3)$ or equivalently if and only if
\begin{gather*}
d=5\text{ and } n\ge 4,\quad d\in\{7,9,11\}\text{ and }n\ge 3,
\quad\text{or}\quad d\ge 13\text{ and }n\ge 2.
\end{gather*}
On the other hand, we compute as in the table
\begin{gather*}
  \begin{array}{|c|ccccc|}
  \hline
 (d,n)&(5,2)&(5,3)&(7,2)&(9,2)&(11,2)\\ 
  \hline
 \dim q_{d,k}(C_{d,n})&d-5=0&d-4=1&d-6=1&d-3=6&d-3=8\\
  \hline
 \end{array}\,.
\end{gather*}
\end{inparaenum}
\end{mainth}

When the odd $d$ equals $3$, we establish a more precise description of
$q_{3,k}(C_{3,n})$ than in the above other $5$ exceptional cases.

\begin{mainth}\label{th:cubic}
For every $n\ge 2$, we have the equality $\dim q_{3,k}(C_{3,n})=0$.
More precisely, 
we have not only
$q_{3,k}(\widevec{\xi_g\infty},1_{\{\widevec{\xi_g\infty}\}}+2\cdot1_{\{\widevec{\xi_g0}\}})\neq
q_{3,k}(\widetilde{(z^{-2})}_{\xi_g},1_{\{\widevec{\xi_g\infty}\}})$ 
but also
\begin{gather}
 q_{3,k}(C_{3,n})=
 \begin{cases}
 \bigl\{q_{3,k}\bigl(\widevec{\xi_g\infty},1_{\{\widevec{\xi_g\infty}\}}+2\cdot1_{\{\widevec{\xi_g0}\}}\bigr)\bigr\}&\text{if }2\le n\le 4,\\
 \bigl\{q_{3,k}\bigl(\widevec{\xi_g\infty},1_{\{\widevec{\xi_g\infty}\}}+2\cdot1_{\{\widevec{\xi_g0}\}}\bigr),q_{3,k}\bigl(\widetilde{(z^{-2})}_{\xi_g},1_{\{\widevec{\xi_g\infty}\}}\bigr)\bigr\}&\text{if }n\ge 5.
 \end{cases}\label{eq:cubic}
\end{gather}
\end{mainth}

For every $n\ge 2$, the above iterate map $\Phi^{\mathrm{IRD}}_{d,n,k}$ descends to the $n$th iterate map
\begin{gather*}
 \Phi_{d,n,k}:\overline{\rM}_{d,k}\dashrightarrow\overline{\rM}_{d^n,k}
\end{gather*}
between the compactified dynamical moduli spaces.
In this compactified dynamical moduli setting, 
one of our principal results is the following characterization 
of the indeterminacy locus $I(\Phi_{d,n,k})$ of $\Phi_{d,n,k}$.
When $k=\bC$, 
Theorem \ref{th:moduli} is due to Kiwi--Nie \cite[Theorem A]{KN23} and 
answered a question of DeMarco \cite[p. 330, Question]{DeMarcoJAMS}, 
based on a sophisticated construction of a degenerating family of rational functions.

\begin{mainth}\label{th:moduli}
For every $d\ge 2$ and $n\ge 2$,
\begin{gather*}
 I(\Phi_{d,n,k})=
 q_{d,k}\Bigl(\bigl\{(\tilde{\phi}_{\xi_g},\dep_{\,\cdot}\tilde{\phi}_{\xi_g})\in I(d):\xi_g\in\MHRL(\phi)\bigr\}\cup C_{d,n}\Bigr).
\end{gather*}
\end{mainth}

We postpone the proof of Theorem \ref{th:moduli},
the argument of which is much simpler than Kiwi-Nie's thanks to Theorem \ref{th:main}(\ref{item:moduli}),
to that of Theorem \ref{th:coeff} in Appendix.
Our final principal result is the following version of 
Theorem \ref{th:main}.

\begin{maincoro}\label{th:indet}
For every $d\ge 2$ and every $n\ge 2$, we have $\dim I(\Phi_{d,n,k})\le d-2$, the equality in which holds if and only if either
\begin{inparaenum}[(i)]
 \item $d$ is even or 
 \item $d$ is odd and $\ge 5$ and $(4/d)^n<1/(d-3)$.
\end{inparaenum}
\end{maincoro}

We recall background materials in Section \ref{sec:background}, and 
prepare some key computations of intrinsic depths in Section \ref{sec:intrinsic}.
We show Theorem \ref{th:main}(\ref{item:coeff}) in Section \ref{sec:upper}
and Theorem \ref{th:main}(\ref{item:moduli}) in Sections \ref{sec:lower} and \ref{sec:exceptional}.
We show Theorem \ref{th:cubic} and Corollary \ref{th:indet} in Section \ref{sec:cubicindet}.
In the Appendix, we show Theorem \ref{th:coeff} and, in particular, Theorem \ref{th:moduli}.

\section{Background}
\label{sec:background}

Let $K$ be an algebraically closed field that is complete
with respect to a non-trivial non-archimedean absolute value, and
let $k$ be the residue field of $K$ of 
not necessarily characteristic $0$.
We denote by $\sP^1$ the Berkovich projective line over $K$, 
by $\sH^1=\sP^1\setminus\bP^1$, where $\bP^1=\bP^1(K)$, the hyperbolic space of $\sP^1$, by
$\sH^1_{\mathrm{II}}$ the set of all type II points in $\sP^1$, by $\rho$
the hyperbolic metric of $\sH^1$, and by $\xi_g\in\sH^1_{\mathrm{II}}$ the Gauss point in $\sP^1$.
The Berkovich projective line $\sP^1$ is compact and uniquely arcwise connected.
The tangent (or direction) space $T_\xi\sP^1$ of $\sP^1$ at $\xi\in\sP^1$ 
is the set of all components of $\sP^1\setminus\{\xi\}$, and 
each direction $\Vv\in T_\xi\sP^1$ is written as
$U(\Vv)$ when it is more convenient, and also as
$\widevec{\xi\xi'}$ for any $\xi'\in U(\Vv)$.
For every $\xi\in\sH^1_{\mathrm{II}}$, the tangent space $T_\xi\sP^1$
has a natural $k$-projective structure; when $\xi=\xi_g$, this structure
can be even canonical by the identification 
\begin{gather*}
 T_{\xi_g}\sP^1\cong\bP^1(k)(=\Rat_{0,k}),\quad \widevec{\xi_ga}\leftrightarrow\hat{a}.
\end{gather*}
The action of $\PGL(2,K^\circ)$ on $T_{\xi_g}\sP^1$
factors through the surjective reduction map
$\PGL(2,K^\circ)\to\PGL(2,k)$ and is transitive
on ordered triples of distinct directions in $T_{\xi_g}\sP^1$.

\subsection{The minimal hyperbolic resultant locus} 
The hyperbolic resultant function
$\hypRes_\phi:\sH^1\to\bR$ for 
a rational function $\phi\in K(z)$ of degree $d>1$
is defined by
\begin{gather*}
 \hypRes_\phi(\xi):=
 \frac{\rho(\xi,\xi_g)}{2}
 +\frac{\rho(\xi,\phi(\xi)\wedge_{\xi_g}\xi)
       -\int_{\sP^1}\rho(\xi_g,\xi\wedge_{\xi_g}\cdot)
                         \,\rd(\phi^*\delta_{\xi_g})}{d-1}
\end{gather*}
on $(\sH^1,\rho)$. For every $\xi,\xi',\xi''\in\sP^1$, the closed interval $[\xi,\xi']$ is the unique closed arc in $\sP^1$ joining $\xi$ and $\xi'$, and
the triple point $\xi\wedge_{\xi''}\xi'\in\sP^1$ is the unique point of
$[\xi,\xi']\cap[\xi',\xi'']\cap[\xi'',\xi]$.
For every $\xi\in\sH^1$ and every direction $\Vv\in T_\xi\sP^1$, we have
the slope formula \cite[(1.1)]{OkuyamaQuadratic}
\begin{gather}
 \rd_{\Vv}\hypRes_\phi
 =\frac{1}{d-1}\Bigl(-\dep_{\Vv}\tilde{\phi}_\xi
   +\frac{1}{2}\cdot
   \begin{cases}
    d-1&\text{if }\tilde{\phi}_\xi(\Vv)=\Vv,\\
    d+1&\text{otherwise}
   \end{cases}\Bigr),\label{eq:slope}
\end{gather}
where $\rd_{\Vv}$ is the derivative at $\xi$ in direction $\Vv$ with respect to the $\rho$-length parameter. 
Extending $\equiv+\infty$ on $\bP^1$,
the hyperbolic resultant function $\hypRes_\phi$
is convex, (metrically) proper, and piecewise affine, and attains its minimum in $\sH^1$ (\cite[Definition 1.1 and Theorem 2]{OkuyamaGeometry}, where $\hypRes_\phi$ was also called the crucial function for $\phi$
and introduced to give a global, explicit, and hyperbolic geometric interpretation of Rumely's function $\ordRes_\phi$ \cite{RumelyMinimal,RumelyEquivariant}).
 
As in Section \ref{sec:intro}, the minimal hyperbolic resultant locus 
is defined by
\begin{gather*}
 \MHRL(\phi):=\arg\min_{\sH^1}\hypRes_\phi,
\end{gather*}
which is a possibly trivial compact interval in $\sH^1$ having
ends in $\sH^1_{\mathrm{II}}$; when $d$ is even, 
$\MHRL(\phi)$ is always a singleton in $\sH^1_{\mathrm{II}}$ by the slope formula \eqref{eq:slope}.

\subsection{Non-archimedean argument principle}
For every rational function $\phi\in K(z)$ of degree $>0$,
every $\xi\in\sH^1$, every direction
$\Vv\in T_\xi\sP^1$, and every $\xi'\in\sP^1$, 
the non-archimedean argument principle
\begin{gather}
 (\phi^*\delta_{\xi'})\bigl(U(\Vv)\bigr)
 =
 s_{\Vv}(\phi)+
 \begin{cases}
 m_{\Vv}(\phi)&\text{if }\xi'\in U\bigl(\phi_{*,\xi}(\Vv)\bigr),\\
 0&\text{otherwise},
\end{cases}\label{eq:FKN}
\end{gather}
is very useful (see, e.g., \cite[Appendix C]{OkuyamaQuadratic} for a proof), where $m_{\Vv}(\phi)\in\{1,\ldots,\deg\phi\}$ and
$s_{\Vv}(\phi)\in\{0,\ldots,\deg\phi-\deg_{\xi}(\phi)\}$ are,
respectively, the directional local degree and the surplus local
degree of $\phi$ at $\Vv$; writing $\Vv=\widevec{\xi\xi'}$, we have
the restriction $\phi:([\xi,\xi'],\rho)\to([\phi(\xi),\phi(\xi')],\rho)$, 
which is a homeomorphism expanding the metric $\rho$ in $m_{\Vv}(\phi)$
if $\xi'$ is close enough to $\xi$, and
\begin{gather}
 \sum_{\Vv\in T_\xi\sP^1:\,\phi_*\Vv=\Vw}m_{\Vv}(\phi)=\deg_\xi(\phi)
\label{eq:directional}
 \end{gather}
for every $\Vw\in T_{\phi(\xi)}\sP^1$.
The following composition formula, which is a consequence of
the non-archimedean argument principle, is also very useful:
for any $\phi,\chi\in K(z)$ of degree $>0$, every
$\xi\in\sH^1_{\mathrm{II}}$ such that $\phi(\xi)=\xi$,
and every $\Vv\in T_\xi\sP^1$, 
\begin{gather}
 \dep_{\Vv}\widetilde{(\chi\circ\phi)}_\xi
 =(\deg\chi)\cdot\dep_{\Vv}\tilde{\phi}_\xi
 +m_{\Vv}(\phi)\cdot\dep_{\phi_{*,\xi}(\Vv)}\tilde{\chi}_\xi;\label{eq:comp}
\end{gather}
for more details, see \cite[Fact 2.1]{OkuyamaQuadratic}.

\subsection{Algebraic stratification of $\IRD_{d,k}(\xi_g)$
and the several auxiliary varieties}\label{sec:strat}

For any integers $d>0$ and $0\le e\le d$, 
we have the $\PGL(2,K^\circ)$-equivariant injective morphism
\begin{gather}
\begin{gathered}
  \iota_e:\bP\bigl(\Sym^{d-e}((k^2)^\vee)\bigr)\times\Rat_{e,k}\to\IRD_{d,k}(\xi_g),\\
 ([H],\tilde{\phi}_{\xi_g})\mapsto \Bigl(\tilde{\phi}_{\xi_g},
 \sum_{\hat{a}\in\bP^1(k)}(\ord_{\hat{a}}[H=0])1_{\{\widevec{\xi_ga}\}}\Bigr);
\end{gathered}
\label{eq:coeff}
\end{gather}
when $e=0$, $\iota_0$ is the Segre embedding
(\cite[I, Exercise 2.14]{Hartshorne}), and
even when $e>0$,
$\iota_e$ is a locally closed immersion. We dispense with this fact below.
These injective morphisms are very useful in computing the loci related to indeterminacies of the iterate maps.

For any integers $1\le m\le e\le d$, we prepare the auxiliary varieties
\begin{align*}
 \cR_{e,m} &:=\bigl\{(\Vv,\tilde{\phi}_{\xi_g})\in T_{\xi_g}\sP^1\times\Rat_{e,k}: m_{\Vv}(\phi)=m\bigr\}\quad\text{and}\\
 \cR_{e,m}'&:=\Bigl\{\bigl(\Vv,(\tilde{\phi}_{\xi_g},\dep_{\,\cdot}\tilde{\phi}_{\xi_g})\bigr)\in T_{\xi_g}\sP^1\times\IRD_{d,k}(\xi_g):(\Vv,\tilde{\phi}_{\xi_g})\in\cR_{e,m}\text{ and }\dep_{\Vv}\tilde{\phi}_{\xi_g}>\frac{d-m}{2}\Bigr\}
\end{align*}
over $k$, the latter of which is $\emptyset$ if $m\le 2e-d$; for, if there is 
$(\Vv,(\tilde{\phi}_{\xi_g},\dep_{\,\cdot}\tilde{\phi}_{\xi_g}))\in\cR_{e,m}'$, then 
$(d-m)/2<\dep_{\Vv}\tilde{\phi}_{\xi_g}\le d-e$. 

For any integers $1\le e\le d$, also set
\begin{gather*}
 \cP_e:=\bigl\{(\Vv,\tilde{\phi}_{\xi_g})\in T_{\xi_g}\sP^1\times\Rat_{e,k}:
 \tilde{\phi}_{\xi_g}(\Vv)\ne\Vv
 \text{ and }(\tilde{\phi}_{\xi_g})^2(\Vv)=\Vv\bigr\},
\end{gather*}
and for every odd integer $d\ge 5$ and
every integer $1\le m\le (d+3)/2$, also set
\begin{gather*}
 \cR_{(d+3)/2,m}''
 :=\bigl\{(\Vv,\tilde{\phi}_{\xi_g})\in\cR_{(d+3)/2,m}:
 \tilde{\phi}_{\xi_g}(\Vv)=\Vv\bigr\}.
\end{gather*}

\section{Key computations of intrinsic depths}
\label{sec:intrinsic}

We retain the setting in Section \ref{sec:background}. Pick $\phi\in\Rat_{d,K}$, $d\ge 2$. We show two lemmas. In the first, we extend the argument in \cite[Section 3]{OkuyamaQuadratic}, and note that
the additional assumption when $\phi(\xi)\ne\xi$ is equivalent to
$(\tilde{\phi}_{\xi_g},\dep_{\,\cdot}\tilde{\phi}_{\xi_g})\not\in I(d)$
if $\xi=\xi_g$.

\begin{lemma}\label{th:depth}
Let $n\ge 2$ and $\xi\in\sH^1_{\mathrm{II}}\setminus\MHRL(\phi^n)$,
and assume that $\dep_{\widevec{\xi\phi(\xi)}}\tilde{\phi}_\xi=0$
when $\phi(\xi)\ne\xi$. Then there is a unique direction 
$\Vv=\Vv_\phi\in T_\xi\sP^1$ such that $\dep_{\Vv}\widetilde{(\phi^n)}_\xi\ge d^n/2$, and then
\begin{gather*}
 2\dep_{\Vv}\tilde{\phi}_\xi
 >d-
 \begin{cases}
 m_{\Vv}(\phi)&\text{if }\phi(\xi)=\xi\\
 0&\text{if }\phi(\xi)\ne\xi
 \end{cases}\ge 0;
\end{gather*}
if in addition 
\begin{gather*}
 \dep_{\Vv_\phi}\tilde{\phi}_\xi\le\frac{d}{2}
\text{ (so in particular }\phi(\xi)=\xi)
\quad\text{and}\quad\dep_{\,\cdot}\tilde{\phi}_{\xi}=\bigl(d-\deg_{\xi}(\phi)\bigr)1_{\Vv_\phi}\text{ on }T_\xi\sP^1,
\end{gather*}
then $\Vv_\phi$ is periodic 
under $\tilde{\phi}_{\xi}$ having a period $p$ in
$\{1,\ldots,n\}$.
\end{lemma}

\begin{proof}
Fix $n\ge 2$ and
$\xi\in\sH^1_{\mathrm{II}}$.
Suppose first that $\phi(\xi)=\xi$.
If $\xi\not\in\MHRL(\phi^n)$, then there is a direction
$\Vv\in T_\xi\sP^1$ such that $\rd_{\Vv}\hypRes_{\phi^n}<0$,
so that 
\begin{gather*}
 \Bigl\{\Vw\in T_\xi\sP^1:\dep_{\Vw}\widetilde{(\phi^n)}_\xi\ge\frac{d^n}{2}\Bigr\}=\{\Vv\}
\end{gather*}
by the slope formula \eqref{eq:slope} applied to $\phi^n$ and
by the inequality
$\sum_{\Vw\in T_\xi\sP^1}\dep_{\Vw}\widetilde{(\phi^n)}_\xi
=d^n-((\phi^n)^*\delta_\xi)(\{\xi\})<d^n$.
Set $i_0:=\min\{i\in\{1,\ldots,n\}:
\dep_{\Vw}\widetilde{(\phi^i)}_\xi\ge d^i/2
\text{ for some }\Vw\in T_\xi\sP^1\}$. Then $\Vv$ is also the unique direction
$\Vw\in T_\xi\sP^1$ such that
$\dep_{\Vw}\widetilde{(\phi^{i_0})}_\xi\ge d^{i_0}/2$;
indeed, for every $\Vw\in T_\xi\sP^1$,
by the composition formula \eqref{eq:comp} applied to
$\phi^n=\phi^{n-i_0}\circ\phi^{i_0}$, we have
$\dep_{\Vw}\widetilde{(\phi^n)}_\xi
 \ge d^{n-i_0}\cdot\dep_{\Vw}\widetilde{(\phi^{i_0})}_\xi$,
and hence $\Vw=\Vv$ whenever $\dep_{\Vw}\widetilde{(\phi^{i_0})}_\xi\ge d^{i_0}/2$.
Hence by the composition formula \eqref{eq:comp} applied to
$\phi^{i_0}=\phi^{i_0-1}\circ\phi$, we have
\begin{multline*}
 \frac{d^{i_0}}{2}\le \dep_{\Vv}\widetilde{(\phi^{i_0})}_\xi
 =\bigl(\deg(\phi^{i_0-1})\bigr)\cdot\dep_{\Vv}\tilde{\phi}_\xi+
 m_{\Vv}(\phi)\cdot\dep_{\tilde{\phi}_\xi(\Vv)}\widetilde{(\phi^{i_0-1})}_\xi\\
 \begin{cases}
 \displaystyle<d^{i_0-1}\cdot\Bigl(\dep_{\Vv}\tilde{\phi}_\xi+\frac{m_{\Vv}(\phi)}{2}\Bigr)
 &\text{if }i_0>1,\\
 \displaystyle=\dep_{\Vv}\tilde{\phi}_\xi+0
 <\dep_{\Vv}\tilde{\phi}_\xi+\frac{m_{\Vv}(\phi)}{2} &\text{if }i_0=1,
 \end{cases}
\end{multline*}
so the desired lower estimate also holds.
If in addition $\dep_{\,\cdot}\tilde{\phi}_{\xi}=\bigl(d-\deg_{\xi}(\phi)\bigr)1_{\Vv}$, $\dep_{\Vv}\tilde{\phi}_\xi\le d/2$,
and $(\tilde{\phi}_\xi)^i(\Vv)\neq\Vv$ for every $i\in\{1,\ldots,n-1\}$,
then by the composition formula \eqref{eq:comp} applied to
$\phi^{i+1}=\phi\circ\phi^i$, we compute
$\dep_{\Vv}\widetilde{(\phi^{i+1})}_\xi
 =d\cdot\dep_{\Vv}\widetilde{(\phi^i)}_\xi
   +m_{\Vv}(\phi^i)\cdot
    \dep_{(\tilde{\phi}_\xi)^i(\Vv)}\tilde{\phi}_\xi
=d\cdot\dep_{\Vv}\widetilde{(\phi^i)}_\xi+0$
for every $i\in\{1,\ldots,n-1\}$, and hence
$\dep_{\Vv}\widetilde{(\phi^n)}_\xi
=d^{n-1}\cdot\dep_{\Vv}\tilde{\phi}_\xi\le d^n/2$.
Unless $(\tilde{\phi}_\xi)^n(\Vv)=\Vv$,
by the slope formula \eqref{eq:slope} applied to $\phi^n$, we must also have
\begin{gather*}
\rd_{\Vv}\hypRes_{\phi^n}
 \ge\frac{1}{d^n-1}\Bigl(-\frac{d^n}{2}+\frac{d^n+1}{2}\Bigr)>0,
\end{gather*}
which is impossible.

Suppose next that $\phi(\xi)\ne\xi$ and that
$\dep_{\widevec{\xi\phi(\xi)}}\tilde{\phi}_\xi=0$,
and set $A=A_\phi:=\{\Vw\in T_\xi\sP^1:
\dep_{\Vw}\tilde{\phi}_\xi>0\}$.
By the non-archimedean argument principle \eqref{eq:FKN}, we have
$A\subset T_\xi\sP^1\setminus\{\widevec{\xi\phi(\xi)}\}$ and
$\phi^{-1}\bigl(U(\widevec{\phi(\xi)\xi})\bigr)
 \subset\bigsqcup_{\Vw\in A}U(\Vw)
 \subset U(\widevec{\phi(\xi)\xi})$.
Hence for every $j\ge 1$, we inductively have
both $\phi^j(\xi)\in\sP^1\setminus U(\widevec{\phi(\xi)\xi})
\subset U(\widevec{\xi\phi(\xi)})$ and 
$\supp((\phi^{j-1})^*\delta_\xi)\subset\{\xi\}\sqcup\bigsqcup_{\Vw\in A}U(\Vw)\subset U(\widevec{\phi(\xi)\xi})$,
and also compute
\begin{align*}
 \dep_{\Vu}\widetilde{(\phi^j)}_\xi
 &=\int_{\sP^1}\bigl(s_{\Vu}(\phi)
   +m_{\Vu}(\phi)\cdot 1_{U(\phi_{*,\xi}(\Vu))}(\cdot)\bigr)
   \bigl((\phi^{j-1})^*\delta_\xi\bigr)(\cdot)\\
 &=d^{j-1}\cdot s_{\Vu}(\phi)
   +\begin{cases}
    m_{\Vu}(\phi)\cdot\bigl((\phi^{j-1})^*\delta_\xi\bigr)
      \bigl(U(\widevec{\phi(\xi)\xi})\bigr)
      &\text{if }\phi_{*,\xi}(\Vu)=\widevec{\phi(\xi)\xi},\\
    0&\text{otherwise}
   \end{cases}\\
 &=d^{j-1}\cdot\dep_{\Vu}\tilde{\phi}_\xi\quad\text{for every }\Vu\in T_\xi\sP^1,
\end{align*}
so that for every $j\ge 1$,
we have $\widetilde{(\phi^j)}_\xi\equiv\widevec{\xi\phi(\xi)}$ and
$\dep_{\widevec{\xi\phi(\xi)}}\widetilde{(\phi^j)}_\xi=0$.
Unless $\xi\in\MHRL(\phi^n)$, also by the slope formula \eqref{eq:slope} applied to $\phi^n$,
there is a direction
$\Vv\in(T_\xi\sP^1)\setminus\{\widevec{\xi\phi(\xi)}\}$ 
such that $\rd_{\Vv}\hypRes_{\phi^n}<0$, and then by the above computation
and the slope formula \eqref{eq:slope} again,
we still have 
\begin{gather*}
 d^{n-1}\cdot\dep_{\Vv}\tilde{\phi}_\xi
 =\dep_{\Vv}\widetilde{(\phi^n)}_\xi
 =\frac{d^n+1}{2}-(d^n-1)\rd_{\Vv}\hypRes_{\phi^n}
 >\frac{d^n+1}{2},
\end{gather*}
so that $\dep_{\Vv}\tilde{\phi}_\xi>d/2$ and in turn
$\dep_{\Vv}\widetilde{(\phi^n)}_\xi>d^n/2$ so
$\sum_{\Vw\in(T_\xi\sP^1)\setminus\{\Vv\}}
\dep_{\Vw}\widetilde{(\phi^n)}_\xi
<d^n-((\phi^n)^*\delta_\xi)(\{\xi\})-d^n/2\le d^n/2$.
Hence $\Vv$ is the unique direction $\Vw\in T_\xi\sP^1$ such that
$\dep_{\Vw}\widetilde{(\phi^n)}_\xi\ge d^n/2$.
\end{proof}

We next consider the case where $\MHRL(\phi)$ is a
non-trivial closed interval.

\begin{lemma}\label{th:nearby}
Let $\xi\in\sH^1_{\mathrm{II}}$
and $\Vv\in T_\xi\sP^1$ satisfy either
\begin{enumerate}[(i)]
\item\label{item:fixed} 
$\phi(\xi)=\xi$, $\tilde{\phi}_\xi(\Vv)=\Vv$,
      $\dep_{\Vv}\tilde{\phi}_\xi=(d-1)/2$, and $m_{\Vv}(\phi)>1$ or
\item\label{item:notfixed} 
$\tilde{\phi}_\xi(\Vv)\neq\Vv$ and
      $\dep_{\Vv}\tilde{\phi}_\xi=(d+1)/2$.
\end{enumerate}
Then there is $\xi_0\in U(\Vv)$ such that
no $\xi'\in(\xi,\xi_0]\cap\sH^1_{\mathrm{II}}$
is fixed by $\phi$ and, for every $\xi'\in(\xi,\xi_0]$,
there is
$\Vw\in T_{\xi'}\sP^1\setminus\{\widevec{\xi'\phi(\xi')}\}$ such that
for every $\Vu\in T_{\xi'}\sP^1$,
\begin{gather*}
 \dep_{\Vu}\tilde{\phi}_{\xi'}=
 \frac{1}{2}\begin{cases}
 d-1&\text{if }\Vu=\widevec{\xi'\phi(\xi')},\\
 d+1&\text{if }\Vu=\Vw,\\
 0&\text{otherwise}.
 \end{cases}
\end{gather*}
\end{lemma}

\begin{proof}
Pick any $\xi\in\sH^1_{\mathrm{II}}$ and any $\Vv\in T_\xi\sP^1$
such that either (\ref{item:fixed}) or (\ref{item:notfixed}) holds. Then there is $\xi_0\in U(\Vv)$ so close to $\xi$ that
when (\ref{item:fixed}) holds, the restriction 
$\phi:([\xi,\xi_0],\rho)\to([\xi,\phi(\xi_0)],\rho)$ 
is a strictly expanding homeomorphism so
$\xi'\in(\xi,\phi(\xi'))$ for every $\xi'\in(\xi,\xi_0]$,
and when (\ref{item:notfixed}) holds, 
for every $\xi'\in(\xi,\xi_0]$, $\widevec{\xi\phi(\xi')}\neq\Vv$
so $\xi\in(\xi',\phi(\xi'))$, and that
$(\phi^*\delta_\xi)(U(\Vv)\cap U(\widevec{\xi_0\xi}))=0$
diminishing $[\xi,\xi_0)$ if necessary.
Then for every $\xi'\in(\xi,\xi_0)$ and
every $\Vu\in T_{\xi'}\sP^1$,
by the non-archimedean argument principle \eqref{eq:FKN}
applied to $\Vu$ and $\delta_{\xi'},\delta_\xi$, we compute
\begin{align*}
 \dep_{\Vu}\tilde{\phi}_{\xi'}
 &=(\phi^*\delta_{\xi'})\bigl(U(\Vu)\bigr)
 =s_{\Vu}(\phi)+
 \begin{cases}
 m_{\Vu}(\phi)&\text{if }\phi_{*,\xi'}(\Vu)=\widevec{\phi(\xi')\xi'}=\widevec{\phi(\xi')\xi},\\
 0&\text{otherwise}
 \end{cases}\\
 &=(\phi^*\delta_\xi)\bigl(U(\Vu)\bigr)
 =\begin{cases}
 \bigl(d-(\phi^*\delta_\xi)(U(\Vv))\bigr)+0
 =d-\dep_{\Vv}\tilde{\phi}_\xi
 &\text{if }\Vu=\widevec{\xi'\xi},\\
 (\phi^*\delta_\xi)(U(\Vv))-0
 =\dep_{\Vv}\tilde{\phi}_\xi
 &\text{if }\Vu=\widevec{\xi'\xi_0},\\
 0&\text{otherwise},
 \end{cases}
\end{align*}
and we are done by setting $\Vw=\widevec{\xi'\xi}$ when (\ref{item:fixed}) holds
and $\Vw=\widevec{\xi'\xi_0}$ when (\ref{item:notfixed}) holds.
\end{proof}

\section{Proof of Theorem \ref{th:main}(\ref{item:coeff})}\label{sec:upper}

Let $d>1$ and $K$ be as in Section \ref{sec:intro}.

\begin{lemma}\label{th:ramif}
For every $e\in\{1,\ldots,d\}$ and every $m\in\{1,\ldots,e\}$,
$\cR_{e,m}$ is irreducible and is of dimension $2e+3-m$, and
for every $p\ge 1$, the subvariety  
$\bigl\{(\Vv,\tilde{\phi}_{\xi_g})\in\cR_{e,m}:(\tilde{\phi}_{\xi_g})^p(\Vv)=\Vv\bigr\}$ in $\cR_{e,m}$ is of dimension $<2e+3-m$.
\end{lemma}

\begin{proof}
Pick $1\le m\le e\le d$.
Let us see that the morphism $\cR_{e,m}\to T_{\xi_g}\sP^1\times T_{\xi_g}\sP^1$, $(\Vv,\tilde{\phi}_{\xi_g})\mapsto(\Vv,\tilde{\phi}_{\xi_g}(\Vv))$, is locally trivial.
Over every $(\Vv,\Vw)\in T_{\xi_g}\sP^1\times T_{\xi_g}\sP^1$, the fiber of this morphism is defined by the equation $\tilde{\phi}_{\xi_g}(\Vv)=\Vw$, and the morphism
\begin{gather*}
 \Bigl\{(\hat{a},\hat{b})\in\bA_k^{e+1}\times\bA_k^{e-m}:
 \hat{a}_e\ne 0\text{ and }
 \gcd\Bigl(\sum_{i=0}^e\hat{a}_i\zeta^i,
 \sum_{i=0}^{e-m-1}\hat{b}_i\zeta^i+\zeta^{e-m}\Bigr)\in k^\times \Bigr\}
\to T_{\xi_g}\sP^1\times\Rat_{e,k},\\
 (\hat{a},\hat{b})=(\hat{a}_0,\ldots,\hat{a}_e,\hat{b}_0,\ldots,\hat{b}_{e-m-1})\mapsto\biggl(\Vv,\widetilde{\Bigl(N^{-1}\circ\frac{\sum_{i=0}^ea_i z^i}{\sum_{i=0}^{e-m-1}b_i z^i+z^{e-m}}\circ M\Bigr)}_{\xi_g}\biggr)
\end{gather*}
is an isomorphism onto this fiber, 
so that the dimension of this fiber is
\begin{gather}
 (e+1)+(e-m)=2e+1-m,\label{eq:fib}
\end{gather}
where we choose $M=M_{\Vv},N=N_{\Vw}\in\PGL(2,K^\circ)$ so that $\tilde{M}_{\xi_g}(\Vv)=\tilde{N}_{\xi_g}(\Vw)=\widevec{\xi_g\infty}$.
Moreover, for every $a\in K^\circ$, 
the tangent map $T_{\xi_g}\sP^1\to T_{\xi_g}\sP^1$ of the projective transformation
\begin{gather*}
 z\mapsto\frac{1}{z-a}\in\PGL(2,K^\circ)
\end{gather*}
of $\bP^1(K^{\circ})$ maps $\widevec{\xi_g a}$ to $\widevec{\xi_g\infty}$,
so we have the desired local triviality by choosing $M,N$ in this form.
Since the base and the fibers are irreducible, $\cR_{e,m}$ is irreducible, and we have
\begin{gather*}
 \dim\cR_{e,m}=(1+1)+(2e+1-m)=2e+3-m.
\end{gather*}
Finally, when $e=m=1$, for every $p\ge 1$, there is $a\in(K^\circ)^\times$ such that $\hat{a}^p\ne1$. Then $(\widevec{\xi_g1},\widetilde{(az)}_{\xi_g})\in\cR_{1,1}$ and $(\widetilde{(az)}_{\xi_g})^p(\widevec{\xi_g1})=\widevec{\xi_g(a^p)}\neq\widevec{\xi_g1}$. 
When $2\le e\le d$ and $1\le m\le e$, we have not only
$(\widevec{\xi_g\infty},\widetilde{(z^{e-m}/(1+z^e-z))}_{\xi_g})\in\cR_{e,m}$ but, for every $p\ge 1$, also $\widetilde{(z^{e-m}/(1+z^e-z))}_{\xi_g}^p(\widevec{\xi_g\infty})\neq\widevec{\xi_g\infty}$. Hence we are done also by the irreducibility of $\cR_{e,m}$.
\end{proof}

\begin{proposition}\label{th:upper}
For every $d\ge 2$ and every $n\ge 2$,
\begin{gather*}
 \dim\bigl\{(\tilde{\phi}_{\xi_g},\dep_{\,\cdot}\tilde{\phi}_{\xi_g})\in\IRD_{d,k}(\xi_g):
 \xi_g\not\in\MHRL(\phi^n)\bigr\}\le d+1,
\end{gather*}
so in particular $\dim C_{d,n}\le d+1$.
\end{proposition}

\begin{proof}
Pick any integers $d,n\ge 2$ and any $\phi\in\Rat_{d,K}$ such that
$\xi_g\notin\MHRL(\phi^n)$.
If $\phi(\xi_g)\ne\xi_g$, then by the Segre embedding $\iota_0$, we have
\begin{multline*}
\dim\bigl\{(\tilde{\phi}_{\xi_g},\dep_{\,\cdot}\tilde{\phi}_{\xi_g})\in\IRD_{d,k}(\xi_g):\phi(\xi_g)\neq\xi_g\bigr\}\\
=\dim\bigl(\bP(\Sym^d((k^2)^\vee))\times\bP^1_k\bigr)=\bigl((d+1)-1\bigr)+1=d+1.
\end{multline*}
Hence suppose that $\phi(\xi_g)=\xi_g$. Then by Lemma \ref{th:depth},
the fiber of the projection $\cR_{\deg_{\xi_g}(\phi),m_{\Vv_{\phi}}(\phi)}'\to\cR_{\deg_{\xi_g}(\phi),m_{\Vv_{\phi}}(\phi)}$
over each $(\Vv,\tilde{\psi}_{\xi_g})\in\cR_{\deg_{\xi_g}(\phi),m_{\Vv_\phi}(\phi)}$ is non-empty and in turn is of dimension 
$<d-\deg_{\xi_g}(\phi)-(d-m_{\Vv_\phi}(\phi))/2$,
so that also using Lemma \ref{th:ramif}, we have
\begin{multline}\label{eq:dimRprime}
 \dim\cR_{\deg_{\xi_g}(\phi),m_{\Vv_\phi}(\phi)}'
 <\dim\cR_{\deg_{\xi_g}(\phi),m_{\Vv_\phi}(\phi)}+\biggl(d-\deg_{\xi_g}(\phi)-\frac{d-m_{\Vv_\phi}(\phi)}{2}\biggr)\\
\le\bigl(2\deg_{\xi_g}(\phi)+3-m_{\Vv_\phi}(\phi)\bigr)
 +\biggl(d-\deg_{\xi_g}(\phi)-\frac{d-m_{\Vv_\phi}(\phi)}{2}\biggr)
 =\deg_{\xi_g}(\phi)+3+\frac{d-m_{\Vv_\phi}(\phi)}{2},\\
\text{which is }\le d+2\text{ if }
m_{\Vv_\phi}(\phi)\ge2\deg_{\xi_g}(\phi)-d+2.
\end{multline}
In the remaining case of $m_{\Vv_\phi}(\phi)=2\deg_{\xi_g}(\phi)-d+1$, 
we have $\dep_{\Vv_\phi}\tilde{\phi}_{\xi_g}\le d/2$
and $\dep_{\,\cdot}\tilde{\phi}_{\xi_g}=(d-\deg_{\xi_g}(\phi))1_{\Vv_\phi}$ on $T_{\xi_g}\sP^1$ since
\begin{gather*}
 d-\deg_{\xi_g}(\phi)-\frac12
 =\frac{d-m_{\Vv_\phi}(\phi)}{2}
 <\dep_{\Vv_\phi}\tilde{\phi}_{\xi_g}
 \le d-\deg_{\xi_g}(\phi)\le\frac{d}{2}.
\end{gather*}
Hence for every $n\ge 2$, also using Lemma \ref{th:depth}, we still have
\begin{gather*}
 \dim\Bigl\{(\tilde{\phi}_{\xi_g},\dep_{\,\cdot}\tilde{\phi}_{\xi_g})
 \in\IRD_{d,k}(\xi_g):
 \xi_g\notin\MHRL(\phi^n)\text{ and }\phi(\xi_g)=\xi_g\Bigr\}\\
\le
\max\biggl\{(d+2)-1,\dim\bigcup_{p=1}^n\bigcup_{1\le 2e-d+1\le e\le d}\Bigl\{\bigl(\Vv,(\tilde{\phi}_{\xi_g},\dep_{\,\cdot}\tilde{\phi}_{\xi_g})\bigr)\in\cR_{e,2e-d+1}':
(\tilde{\phi}_{\xi_g})^p(\Vv)=\Vv
\Bigr\}\biggr\}\\
 \le\max\Bigl\{d+1,\max_{1\le 2e-d+1\le e\le d}\bigl(2e+3-(2e-d+1)\bigr)-1\Bigr\}=d+1,
\end{gather*}
and we are done.
\end{proof}

\begin{proposition}\label{th:evenlower}
For every even $d\ge 2$ and every $n\ge 2$, we have $\dim C_{d,n}\ge d+1$.
\end{proposition}

\begin{proof}
First, let us see the local triviality of the morphism
\begin{gather*}
 \cP_{d/2}\to(T_{\xi_g}\sP^1\times T_{\xi_g}\sP^1)\setminus\Delta_{T_{\xi_g}\sP^1},
 \quad(\Vv,\tilde{\phi}_{\xi_g})\mapsto(\Vv,\tilde{\phi}_{\xi_g}(\Vv)),
\end{gather*}
which yields 
$\dim\cP_{d/2}=(1+1)+((d/2)+(d/2-1))=d+1$.
Over every $(\Vv,\Vw)\in(T_{\xi_g}\sP^1\times T_{\xi_g}\sP^1)\setminus\Delta_{T_{\xi_g}\sP^1}$,
the fiber of this morphism is defined by the equations
$\tilde{\phi}_{\xi_g}(\Vv)=\Vw$ and $\tilde{\phi}_{\xi_g}(\Vw)=\Vv$,
and the morphism
\begin{gather*}
 \Bigl\{(\hat{a},\hat{b})\in\bA_k^{d/2}\times\bA_k^{d/2-1}:
 \sum_{j=0}^{d/2-1}\hat{a}_j\zeta^j\text{ and }
 \zeta^{d/2}+\sum_{j=1}^{d/2-1}\hat{b}_j\zeta^j\text{ are coprime}
 \Bigr\}\to T_{\xi_g}\sP^1\times\Rat_{d/2,k},\\
 (\hat{a},\hat{b})=(\hat{a}_0,\ldots,\hat{a}_{d/2-1},\hat{b}_1,\ldots,\hat{b}_{d/2-1})
 \mapsto\left(\Vv,\widetilde{\Bigl(M^{-1}\circ
 \frac{\sum_{j=0}^{d/2-1}a_j z^j}
 {z^{d/2}+\sum_{j=1}^{d/2-1}b_j z^j}\circ M\Bigr)}_{\xi_g}\right)
\end{gather*}
is an isomorphism onto this fiber, where we choose
$M=M_{(\Vv,\Vw)}\in\PGL(2,K^\circ)$ so that
$\tilde{M}_{\xi_g}(\Vv)=\widevec{\xi_g\infty}$ and
$\tilde{M}_{\xi_g}(\Vw)=\widevec{\xi_g0}$;
the source of the morphism is nonempty since $1$ and $\zeta^{d/2}$ are coprime. 
Moreover, for every $a,b\in K^\circ$ such that $\hat{a}\ne\hat{b}$,
the tangent map of the projective transformation
\begin{gather*}
 z\mapsto\frac{z-b}{z-a}\in\PGL(2,K^\circ)
\end{gather*}
sends $\widevec{\xi_ga}$ and $\widevec{\xi_gb}$ to
$\widevec{\xi_g\infty}$ and $\widevec{\xi_g0}$, respectively,
so we have the desired local triviality by choosing $M$ in such a form. 

Next, we have the injective morphism
\begin{gather*}
 \kappa_{d/2}:\cP_{d/2}\to\IRD_{d,k}(\xi_g),\quad
 (\Vv,\tilde{\phi}_{\xi_g})\mapsto\bigl(\tilde{\phi}_{\xi_g},
\frac{d}{2}\cdot 1_{\{\Vv\}}\bigr)
\end{gather*}
using $\iota_{d/2}$ in \eqref{eq:coeff},
and let us see that $\kappa_{d/2}(\cP_{d/2})\subset C_{d,n}$
for every $n\ge 2$,
which will complete the proof; indeed, $\kappa_{d/2}(\cP_{d/2})\subset\IRD_{d,k}(\xi_g)\setminus I(d)$ since $\deg_{\xi_g}(\phi)=d-d/2=d/2>0$.
Pick any $\phi\in\Rat_{d,K}$ such that
$(\tilde{\phi}_{\xi_g},\dep_{\,\cdot}\tilde{\phi}_{\xi_g})\in \kappa_{d/2}(\cP_{d/2})$. Then $\{\xi_g\}=\MHRL(\phi)$ by the slope formula \eqref{eq:slope} applied to $\phi$ at $\xi_g$ for each $\Vw\in T_{\xi_g}\sP^1$.
Moreover, for $n=2$, recalling $\tilde{\phi}_{\xi_g}(\Vv)\neq\Vv$,
by the composition formula \eqref{eq:comp} applied to $\phi^2=\phi\circ\phi$, we have
\begin{gather*}
 \dep_{\Vv}\widetilde{(\phi^2)}_{\xi_g}
 =(\deg\phi)\cdot\dep_{\Vv}\tilde{\phi}_{\xi_g}
 +m_{\Vv}(\phi)\cdot\dep_{\tilde{\phi}_{\xi_g}(\Vv)}\tilde{\phi}_{\xi_g}
 =d\cdot\frac{d}{2}+0=\frac{d^2}{2},
\end{gather*}
so that also recalling $(\tilde{\phi}_{\xi_g})^2(\Vv)=\Vv$,
by the slope formula \eqref{eq:slope} applied to $\phi^2$,
\begin{gather*}
 \rd_{\Vv}\hypRes_{\phi^2}=\frac{1}{d^2-1}\Bigl(-\frac{d^2}{2}+\frac{d^2-1}{2}\Bigr)<0,
\end{gather*}
and similarly, for every $n\ge 3$, by the composition formula \eqref{eq:comp} applied to
$\phi^n=\phi^{n-2}\circ\phi^2$ and $\phi^{n-2}=\phi^{n-3}\circ\phi$, we have
\begin{multline*}
 \dep_{\Vv}\widetilde{(\phi^n)}_{\xi_g}
 =\deg(\phi^{n-2})\cdot\dep_{\Vv}\widetilde{(\phi^2)}_{\xi_g}
 +m_{\Vv}(\phi^2)\cdot\dep_{\Vv}\widetilde{(\phi^{n-2})}_{\xi_g}\\
 \ge d^{n-2}\cdot\frac{d^2}{2}+1\cdot\Bigl(\deg(\phi^{n-3})\cdot\dep_{\Vv}\widetilde{\phi}_{\xi_g}+0\Bigr)=\frac{d^n}{2}+d^{n-3}\cdot\frac{d}{2},
\end{multline*}
so that by the slope formula \eqref{eq:slope} applied to $\phi^n$, 
\begin{gather*}
 \rd_{\Vv}\hypRes_{\phi^n}
 \le\frac{1}{d^n-1}\Bigl(-\Bigl(\frac{d^n}{2}+d^{n-3}\cdot\frac{d}{2}\Bigr)
 +\frac{d^n+1}{2}\Bigr)<0.
\end{gather*}
Hence $\xi_g\not\in\MHRL(\phi^n)$ for every $n\ge 2$, and we are done.
\end{proof}

The proof of the following is similar to but simpler than that in the case that $d$ is even.

\begin{proposition}\label{th:oddlower}
For every odd $d\ge 3$ and every $n\ge 2$, we have $\dim C_{d,n}\ge d+1$.
\end{proposition}

\begin{proof}
First, the pairs $(\Vv,\tilde{\phi}_{\xi_g})$ such that
$\tilde{\phi}_{\xi_g}(\Vv)\ne\Vv$ form an open subset of
$T_{\xi_g}\sP^1\times\Rat_{(d-1)/2,k}$ containing $\cP_{(d-1)/2}$,
which is nonempty by the proof of Proposition \ref{th:evenlower}.
This open subset is of dimension
$\dim(T_{\xi_g}\sP^1)+\dim(\bP(V_{(d-1)/2,k}))=d+1$.

Next, we have the injective morphism
\begin{gather*}
 \kappa_{(d-1)/2}:T_{\xi_g}\sP^1\times\Rat_{(d-1)/2,k}\to\IRD_{d,k}(\xi_g),\quad (\Vv,\tilde{\phi}_{\xi_g})\mapsto\biggl(\tilde{\phi}_{\xi_g},
 \frac{d+1}{2}1_{\{\Vv\}}\biggr)
\end{gather*}
using $\iota_{(d-1)/2}$ in \eqref{eq:coeff}.
Pick any $(\Vv,\tilde{\phi}_{\xi_g})$ in this open subset and any
$\phi\in\Rat_{d,K}$ such that
$(\tilde{\phi}_{\xi_g},\dep_{\,\cdot}\tilde{\phi}_{\xi_g})
=\kappa_{(d-1)/2}(\Vv,\tilde{\phi}_{\xi_g})$,
which is not in $I(d)$ since $\deg_{\xi_g}(\phi)=d-(d+1)/2=(d-1)/2>0$.
Then $\xi_g\in\MHRL(\phi)$ by the slope formula \eqref{eq:slope}
applied to $\phi$ at $\xi_g$ for each $\Vw\in T_{\xi_g}\sP^1$.
Moreover, for every $n\ge 2$,
by the composition formula \eqref{eq:comp} applied to
$\phi^n=\phi^{n-1}\circ\phi$, we have
\begin{multline*}
 \dep_{\Vv}\widetilde{(\phi^n)}_{\xi_g}
 =\deg(\phi^{n-1})\cdot\dep_{\Vv}\tilde{\phi}_{\xi_g}
 +m_{\Vv}(\phi)\cdot\dep_{\tilde{\phi}_{\xi_g}(\Vv)}\widetilde{(\phi^{n-1})}_{\xi_g}\\
 \ge d^{n-1}\cdot\frac{d+1}{2}+0
 =\frac{d^n+d^{n-1}}{2},
\end{multline*}
so that also using the slope formula \eqref{eq:slope} applied to $\phi^n$, we have
\begin{gather*}
 \rd_{\Vv}\hypRes_{\phi^n}
 \le\frac{1}{d^n-1}\Bigl(-\frac{d^n+d^{n-1}}{2}+\frac{d^n+1}{2}\Bigr)<0.
\end{gather*}
Hence $\xi_g\not\in\MHRL(\phi^n)$, and we are done.
\end{proof}
Now the proof of Theorem \ref{th:main}(\ref{item:coeff})
for every $n\ge 2$ is complete.

\section{Proof of Theorem \ref{th:main}(\ref{item:moduli}):
the non-exceptional cases}\label{sec:lower}

Let $d>1$ and $K$ be as in Section \ref{sec:intro}. Recall the action 
\begin{gather}
 \Bigl(M,\bigl(\tilde{\phi}_{\xi_g},\dep_{\,\cdot}\tilde{\phi}_{\xi_g}\bigr)\Bigr)\mapsto\Bigl(\widetilde{(M\circ\phi\circ M^{-1})}_{\xi_g},
 \dep_{\tilde{M}_{\xi_g}^{-1}}\tilde{\phi}_{\xi_g}\Bigr)\label{eq:stable}
\end{gather}
on $\IRD_{d,k}(\xi_g)$ of $\PGL(2,K^\circ)$; see Subsection \ref{sec:form}.

For every $\phi\in\Rat_{d,K}$,
if $\MHRL(\phi)=\{\xi_g\}$ or equivalently
$(\tilde{\phi}_{\xi_g},\dep_{\,\cdot}\tilde{\phi}_{\xi_g})
\in\IRD_{d,k}^\s(\xi_g)$, then by the definition of a geometric quotient,
\begin{gather}
 \dim q_{d,k}^{-1}\bigl(q_{d,k}(\tilde{\phi}_{\xi_g},
 \dep_{\,\cdot}\tilde{\phi}_{\xi_g})\bigr)
 =\dim\Bigl(\PGL(2,K^\circ)\cdot
 (\tilde{\phi}_{\xi_g},\dep_{\,\cdot}\tilde{\phi}_{\xi_g})\Bigr)
 =\dim\PGL_{2,k}=3.\label{eq:fibers}
\end{gather}
When $d$ is even, $\MHRL(\phi)=\{\xi_g\}$ for every
$(\tilde{\phi}_{\xi_g},\dep_{\,\cdot}\tilde{\phi}_{\xi_g})\in C_{d,n}$
by the slope formula \eqref{eq:slope}, 
so in particular for every $n\ge 2$, the desired equality
\begin{gather*}
 \dim q_{d,k}(C_{d,n})=\dim C_{d,n}-3=(d+1)-3=d-2
\end{gather*}
holds by \eqref{eq:fibers} and Theorem \ref{th:main}(\ref{item:coeff}).

Let us focus on the case that $d$ is odd.
The following concludes the desired upper estimate
of $\dim q_{d,k}(C_{d,n})$ for every $n\ge 2$ since
by Proposition \ref{th:upper} and \eqref{eq:fibers},
$\dim q_{d,k}(\{(\tilde{\phi}_{\xi_g},\dep_{\,\cdot}\tilde{\phi}_{\xi_g})
\in C_{d,n}:\MHRL(\phi)=\{\xi_g\}\})\le (d+1)-3=d-2$.

\begin{lemma}\label{th:stab}
For every odd $d\ge 3$, every $n\ge 2$, and every $\phi\in\Rat_{d,K}$ such that
$(\tilde{\phi}_{\xi_g},\dep_{\,\cdot}\tilde{\phi}_{\xi_g})\in C_{d,n}$
and $\{\xi_g\}\subsetneq\MHRL(\phi)$, we have
\begin{gather*}
 q_{d,k}\bigl(\tilde{\phi}_{\xi_g},\dep_{\,\cdot}\tilde{\phi}_{\xi_g}\bigr)
 =q_{d,k}\Bigl(\widevec{\xi_g\infty},
 \frac{d-1}{2}1_{\{\widevec{\xi_g\infty}\}}
 +\frac{d+1}{2}1_{\{\widevec{\xi_g0}\}}\Bigr).
\end{gather*}
\end{lemma}

\begin{proof}
Pick any $\phi\in\Rat_{d,K}$ such that
$(\tilde{\phi}_{\xi_g},\dep_{\,\cdot}\tilde{\phi}_{\xi_g})\in C_{d,n}$ and
that $\{\xi_g\}\subsetneq\MHRL(\phi)$. We claim that
either
(\ref{item:fixed}) or (\ref{item:notfixed}) in Lemma \ref{th:nearby} 
is the case. Indeed, by the slope formula \eqref{eq:slope}, 
$\rd_{\Vv}\hypRes_{\phi}=0$ for some $\Vv\in T_{\xi_g}\sP^1$, and then
\begin{gather*}
 \dep_{\Vv}\tilde{\phi}_{\xi_g}=-0+\frac{1}{2}
 \begin{cases}
  d+1 &\text{if }\tilde{\phi}_{\xi_g}(\Vv)\neq\Vv,\\
  d-1 &\text{if }\tilde{\phi}_{\xi_g}(\Vv)=\Vv;
 \end{cases}
\end{gather*}
in the latter case, we have not only 
$\phi(\xi_g)=\xi_g$ since otherwise we must have
$\tilde{\phi}_{\xi_g}\equiv\Vv$ and then
$(\tilde{\phi}_{\xi_g},\dep_{\,\cdot}\tilde{\phi}_{\xi_g})\in I(d)$,
which is impossible, but also $\Vv=\Vv_\phi$ by 
the computation
\begin{multline*}
 \sup_{\Vw\in(T_{\xi_g}\sP^1)\setminus\{\Vv\}}\dep_{\Vw}\tilde{\phi}_{\xi_g}
 \le\sum_{\Vu\in(T_{\xi_g}\sP^1)\setminus\{\Vv\}}\dep_{\Vu}\tilde{\phi}_{\xi_g}\\
 =d-\deg_{\xi_g}(\phi)-\frac{d-1}{2}
 =\frac{d-2\deg_{\xi_g}(\phi)+1}{2}
\le\frac{d-\sup_{\Vw\in(T_{\xi_g}\sP^1)\setminus\{\Vv\}}m_{\Vw}(\phi)}{2}
\end{multline*}
and Lemma \ref{th:depth} applied to $\phi$ at $\xi_g$,
and then we also have $m_{\Vv}(\phi)>1$ by 
$(d-1)/2=\dep_{\Vv}\tilde{\phi}_{\xi_g}>(d-m_{\Vv}(\phi))/2$
and Lemma \ref{th:depth} applied to $\phi$ and $\Vv=\Vv_\phi$.

Once the claim is at our disposal,
by Lemma \ref{th:nearby}, there are also
$\xi'\in(\MHRL(\phi)\setminus\{\xi_g\})\cap\sH^1_{\mathrm{II}}$,
$\Vw\in(T_{\xi'}\sP^1)\setminus\{\widevec{\xi'\phi(\xi')}\}$, 
and $M\in\PGL(2,K)$ such that $M(\xi')=\xi_g$,
$M_{*,\xi'}(\widevec{\xi'\phi(\xi')})=\widevec{\xi_g\infty}$, and
$M_{*,\xi'}(\Vw)=\widevec{\xi_g0}$ and that
\begin{align*}
 q_{d,k}\bigl(\tilde{\phi}_{\xi_g},\dep_{\,\cdot}\tilde{\phi}_{\xi_g}\bigr)
 &=q_{d,k}\Bigl(\widetilde{(M\circ\phi\circ M^{-1})}_{\xi_g},\dep_{\,\cdot}\widetilde{(M\circ\phi\circ M^{-1})}_{\xi_g}\Bigr)\\
 &=q_{d,k}\Bigl(\widevec{\xi_g\infty},
 \frac{d-1}{2}1_{\{\widevec{\xi_g\infty}\}}
 +\frac{d+1}{2}1_{\{\widevec{\xi_g0}\}}\Bigr),
\end{align*}
recalling also the final paragraph of Subsection \ref{sec:form}, and we are done.
\end{proof}

When $d$ is odd and $\ge 5$, for every $m\in\{1,\ldots,(d+3)/2\}$,
we have the injective morphism
\begin{gather*}
 \lambda_{(d+3)/2,m}:\cR_{(d+3)/2,m}''
\to\IRD_{d,k}(\xi_g),\quad
 (\Vv,\tilde{\phi}_{\xi_g})\mapsto\Bigl(\tilde{\phi}_{\xi_g},
 \frac{d-3}{2}1_{\{\Vv\}}\Bigr),
\end{gather*}
using $\iota_{(d+3)/2}$ in \eqref{eq:coeff},
the image of
which not only is invariant under the action \eqref{eq:stable}
but, by $\deg_{\xi_g}(\phi)=(d+3)/2>0$, also is contained
in $\IRD_{d,k}(\xi_g)\setminus I(d)$.

The following will complete the proof of
Theorem \ref{th:main}(\ref{item:moduli})
in the non-exceptional cases.

\begin{proposition}[{cf.\ \cite[Corollary 4.16]{KN23} when $d=5$}]
\label{th:poly}
For every odd integer $d\ge 5$ and every $n\ge 2$ such that
$(4/d)^n<1/(d-3)$, we have $\dim q_{d,k}(C_{d,n})\ge d-2$.
\end{proposition}

\begin{proof}
Pick any odd integer $d\ge 5$. 
The fiber over each $\Vv\in T_{\xi_g}\sP^1$
of the projection $\cR_{(d+3)/2,4}''\to T_{\xi_g}\sP^1$
is nothing but the fiber over $(\Vv,\Vv)\in\Delta_{T_{\xi_g}\sP^1}$
of the morphism
$\cR_{(d+3)/2,4}\to T_{\xi_g}\sP^1\times T_{\xi_g}\sP^1$,
$(\Vv,\tilde{\phi}_{\xi_g})\mapsto(\Vv,\tilde{\phi}_{\xi_g}(\Vv))$
in the proof of Lemma \ref{th:ramif},
so in particular is of dimension $2\cdot(d+3)/2+1-4=d$. Hence
\begin{gather*}
 \dim\cR_{(d+3)/2,4}''
 =\dim(T_{\xi_g}\sP^1)+d=d+1.
\end{gather*}
Next, pick any $n\ge 2$ such that $(4/d)^n<1/(d-3)$, and
pick any $(\Vv,\tilde{\phi}_{\xi_g})\in\cR_{(d+3)/2,4}''$
and any $\phi\in\Rat_{d,K}$ such that
$(\tilde{\phi}_{\xi_g},\dep_{\,\cdot}\tilde{\phi}_{\xi_g})=\lambda_{(d+3)/2,4}(\Vv,\tilde{\phi}_{\xi_g})$. Then we have not only 
$\MHRL(\phi)=\{\xi_g\}$ since by the slope formula \eqref{eq:slope}, 
\begin{gather*}
 \rd_{\Vv}\hypRes_{\phi}
 =\frac{1}{d-1}\Bigl(-\frac{d-3}{2}+\frac{d-1}{2}\Bigr)
 =\frac{1}{d-1}>0
 \quad\text{and}\\
 \rd_{\Vw}\hypRes_{\phi}
 \ge\frac{1}{d-1}\Bigl(-0+\frac{d-1}{2}\Bigr)
 =\frac{1}{2}>0\quad\text{for every }
\Vw\in T_{\xi_g}\sP^1\setminus\{\Vv\},
\end{gather*}
but also $\xi_g\not\in\MHRL(\phi^n)$ since
by the composition formula \eqref{eq:comp} applied to
$\phi^n=\phi^{n-1}\circ\phi$, 
also using the assumption on $d$ and $n$, 
we recursively have
\begin{multline*}
 \dep_{\Vv}\widetilde{(\phi^n)}_{\xi_g}
 =\deg(\phi^{n-1})\cdot\dep_{\Vv}\tilde{\phi}_{\xi_g}
 +m_{\Vv}(\phi)\cdot\dep_{\tilde{\phi}_{\xi_g}(\Vv)}\widetilde{(\phi^{n-1})}_{\xi_g}\\
 =d^{n-1}\cdot\frac{d-3}{2}+4\cdot\dep_{\Vv}\widetilde{(\phi^{n-1})}_{\xi_g}
 =\frac{d-3}{2}\cdot\sum_{j=0}^{n-1}d^{n-1-j}\cdot4^j
 =\frac{d-3}{2(d-4)}\cdot(d^n-4^n)\\
 \ge\frac{d-3}{2(d-4)}\cdot
 \Bigl(d^n-\frac{d^n}{d-3}\Bigr)
 =\frac{d-3}{2(d-4)}
 \cdot\frac{(d-3)d^n-d^n}{d-3}
 =\frac{d^n}{2},
\end{multline*}
which with the slope formula
\eqref{eq:slope} applied to $\phi^n$ yields
\begin{gather*}
 \rd_{\Vv}\hypRes_{\phi^n}
 \le\frac{1}{d^n-1}\Bigl(-\frac{d^n}{2}+\frac{d^n-1}{2}\Bigr)
 <0.
\end{gather*}
Hence
$\dim q_{d,k}(C_{d,n})
 \ge\dim\cR_{(d+3)/2,4}''-3=(d+1)-3=d-2$
by \eqref{eq:fibers}, and we are done.
\end{proof}

\section{Proof of Theorem \ref{th:main}(\ref{item:moduli}):
the five exceptional cases}\label{sec:exceptional}

Let $K$ be as in Section \ref{sec:intro}. The following computation
is useful in treating the remaining exceptional cases.

\begin{lemma}\label{th:dimcycle}
 $\dim(\cP_2\cap\cR_{2,2})=4$.
\end{lemma}

\begin{proof}
 Over every $(\Vv,\Vw)\in(T_{\xi_g}\sP^1\times T_{\xi_g}\sP^1)\setminus\Delta_{T_{\xi_g}\sP^1}$, the fiber of 
 the morphism
 \begin{gather*}
 \cP_2\cap\cR_{2,2}\to(T_{\xi_g}\sP^1\times T_{\xi_g}\sP^1)\setminus\Delta_{T_{\xi_g}\sP^1},
 \quad(\Vv,\tilde{\phi}_{\xi_g})\mapsto(\Vv,\tilde{\phi}_{\xi_g}(\Vv))
 \end{gather*}
 is defined by the equations
 $\tilde{\phi}_{\xi_g}(\Vv)=\Vw$ and $\tilde{\phi}_{\xi_g}(\Vw)=\Vv$.
 Choose $M\in\PGL(2,K^\circ)$ such that
 $\tilde{M}_{\xi_g}(\Vv)=\widevec{\xi_g\infty}$ and
 $\tilde{M}_{\xi_g}(\Vw)=\widevec{\xi_g0}$.
 The morphism
 \begin{gather*}
 (\bA_k^1\setminus\{0\})\times\bA_k^1\to\cP_2\cap\cR_{2,2},
 \quad(\hat{a},\hat{b})\mapsto
 \Bigl(\Vv,\widetilde{\Bigl(M^{-1}\circ\frac{a}{z^2+bz}\circ M\Bigr)}_{\xi_g}\Bigr)
 \end{gather*}
 is an isomorphism onto this fiber, so that
$\dim(\cP_2\cap\cR_{2,2})
 =\dim((T_{\xi_g}\sP^1\times T_{\xi_g}\sP^1)\setminus\Delta_{T_{\xi_g}\sP^1})
 +\dim((\bA_k^1\setminus\{0\})\times\bA_k^1)
 =(1+1)+(1+1)=4$.
\end{proof}

Pick any odd $d\ge 5$ and any $n\ge 2$, and
pick any $\phi\in\Rat_{d,K}$ such that 
\begin{gather*}
 (\tilde{\phi}_{\xi_g},\dep_{\,\cdot}\tilde{\phi}_{\xi_g})\in C_{d,n} 
\quad\text{and}\quad
\MHRL(\phi)=\{\xi_g\}. 
\end{gather*}
Then we have
\begin{gather*}
 \phi(\xi_g)=\xi_g;
\end{gather*}
for, otherwise, by Lemma \ref{th:depth} and the slope formula
\eqref{eq:slope}, we must have
$d/2<\dep_{\Vv_\phi}\tilde{\phi}_{\xi_g}<(d+1)/2$,
and then $\le(d-1)/2$, which is impossible.
Then by Lemma \ref{th:depth} and the slope formula \eqref{eq:slope} again,
we have not only
\begin{gather}
 0\le\frac{d-m_{\Vv_\phi}(\phi)}{2}<
 \dep_{\Vv_\phi}\tilde{\phi}_{\xi_g}<\frac{d+1}{2},
\label{eq:dep}
\end{gather}
so in particular
$(\Vv_\phi,(\tilde{\phi}_{\xi_g},\dep_{\,\cdot}\tilde{\phi}_{\xi_g}))
\in\cR_{\deg_{\xi_g}(\phi),m_{\Vv_\phi}(\phi)}'$, but also
\begin{gather}
 \dep_{\Vv_\phi}\tilde{\phi}_{\xi_g}
<\frac{d-1}{2}
\quad\text{when }\tilde{\phi}_{\xi_g}(\Vv_\phi)=\Vv_\phi.\label{eq:depfix}
\end{gather}
We also note that
$\dep_{\Vv_\phi}\tilde{\phi}_{\xi_g}
 \le\sum_{\Vw\in T_{\xi_g}\sP^1}\dep_{\Vw}\tilde{\phi}_{\xi_g}
 =d-\deg_{\xi_g}(\phi)$.

\begin{inparaenum}[(i)]
\item Suppose that $n=2$.
\begin{inparaenum}[({i}-1)]
\item When $\dep_{\Vv_\phi}\tilde{\phi}_{\xi_g}
\le d-\deg_{\xi_g}(\phi)-2$, using the estimate \eqref{eq:dimRprime} of
$\dim\cR_{\deg_{\xi_g}(\phi),m_{\Vv_\phi}(\phi)}'$ and 
\eqref{eq:dep},
we have
\begin{multline*}
 \dim\cR_{\deg_{\xi_g}(\phi),m_{\Vv_\phi}(\phi)}'
 <\frac{1}{2}\bigl(2\deg_{\xi_g}(\phi)+6+d-m_{\Vv_\phi}(\phi)\bigr)\\
<\frac{1}{2}\Bigl(2\deg_{\xi_g}(\phi)+6+d-
\bigl(d-2\bigl(d-\deg_{\xi_g}(\phi)-2\bigr)\bigr)\Bigr)
 =d+1.
\end{multline*}
\item When $\dep_{\Vv_\phi}\tilde{\phi}_{\xi_g}=d-\deg_{\xi_g}(\phi)-1$,
we have $\tilde{\phi}_{\xi_g}(\Vv_\phi)=\Vv_\phi$;
for, otherwise, we have
$\dep_{\tilde{\phi}_{\xi_g}(\Vv_\phi)}\tilde{\phi}_{\xi_g}
\le\sum_{\Vw\in(T_{\xi_g}\sP^1)\setminus\{\Vv_\phi\}}\dep_{\Vw}\tilde{\phi}_{\xi_g}
=d-\deg_{\xi_g}(\phi)-\dep_{\Vv_\phi}\tilde{\phi}_{\xi_g}=1$,
and in turn by Lemma \ref{th:depth}, the composition formula \eqref{eq:comp} applied to $\phi^2=\phi\circ\phi$, and \eqref{eq:dep},
\begin{multline*}
 \frac{d^2}{2}\le\dep_{\Vv_\phi}\widetilde{(\phi^2)}_{\xi_g}
 =(\deg\phi)\cdot\dep_{\Vv_\phi}\tilde{\phi}_{\xi_g}
 +m_{\Vv_\phi}(\phi)\cdot\dep_{\tilde{\phi}_{\xi_g}(\Vv_\phi)}\tilde{\phi}_{\xi_g}\\
 \le d\bigl(d-\deg_{\xi_g}(\phi)-1\bigr)+\deg_{\xi_g}(\phi)\cdot 1
 =(d-1)(d-\deg_{\xi_g}(\phi))\\
 =(d-1)\bigl(\dep_{\Vv_\phi}\tilde{\phi}_{\xi_g}+1\bigr)
 \le(d-1)\Bigl(\frac{d-1}{2}+1\Bigr)
 =\frac{d^2-1}{2},
\end{multline*}
which is impossible. Then also using the sharper estimate \eqref{eq:depfix},
we have
\begin{multline*}
 \frac{d^2}{2}\le\dep_{\Vv_\phi}\widetilde{(\phi^2)}_{\xi_g}
 =d\bigl(d-\deg_{\xi_g}(\phi)-1\bigr)
 +m_{\Vv_\phi}(\phi)\cdot\dep_{\Vv_\phi}\tilde{\phi}_{\xi_g}\\
 \le d\bigl(d-\deg_{\xi_g}(\phi)-1\bigr)
+\frac{d-3}{2}\cdot m_{\Vv_\phi}(\phi).
\end{multline*}
Hence if $d<12$, then using the sharper estimate \eqref{eq:depfix} again,
we further compute
\begin{multline*}
 m_{\Vv_\phi}(\phi)
 \ge\frac{d(-d+2\deg_{\xi_g}(\phi)+2)}{d-3}
 =-d+2\deg_{\xi_g}(\phi)+2
 +3\frac{d-2\dep_{\Vv_\phi}\tilde{\phi}_{\xi_g}}{d-3}\\
 \ge-d+2\deg_{\xi_g}(\phi)+2+3\frac{d-(d-3)}{d-3}
 >-d+2\deg_{\xi_g}(\phi)+3,
\end{multline*}
and in turn using the estimate \eqref{eq:dimRprime} of
$\dim\cR_{\deg_{\xi_g}(\phi),m_{\Vv_\phi}(\phi)}'$, we have
\begin{multline*}
 \dim\cR_{\deg_{\xi_g}(\phi),m_{\Vv_\phi}(\phi)}'
 <\frac{1}{2}\bigl(2\deg_{\xi_g}(\phi)+6+d-m_{\Vv_\phi}(\phi)\bigr)\\
 \le\frac{1}{2}\bigl(2\deg_{\xi_g}(\phi)+6+d-(-d+2\deg_{\xi_g}(\phi)+4)\bigr)=d+1.
\end{multline*}
\item When $\dep_{\Vv_\phi}\tilde{\phi}_{\xi_g}=d-\deg_{\xi_g}(\phi)$,
that is, $\dep_{\,\cdot}\tilde{\phi}_{\xi_g}
=(d-\deg_{\xi_g}(\phi))1_{\{\Vv_\phi\}}$, we still have
$\tilde{\phi}_{\xi_g}(\Vv_\phi)=\Vv_\phi$;
for, otherwise,
$\dep_{\tilde{\phi}_{\xi_g}(\Vv_\phi)}\tilde{\phi}_{\xi_g}=0$,
and in turn
by Lemma \ref{th:depth}, the composition formula \eqref{eq:comp}
applied to $\phi^2=\phi\circ\phi$, and \eqref{eq:dep}, we must have
\begin{gather*}
 \frac{d^2}{2}\le\dep_{\Vv_\phi}\widetilde{(\phi^2)}_{\xi_g}
 =(\deg\phi)\cdot\dep_{\Vv_\phi}\tilde{\phi}_{\xi_g}
+m_{\Vv_\phi}(\phi)\cdot\dep_{\tilde{\phi}_{\xi_g}(\Vv_\phi)}\tilde{\phi}_{\xi_g}
 \le d\cdot\frac{d-1}{2},
\end{gather*}
which is impossible.
Then also using the sharper estimate \eqref{eq:depfix}, we have
\begin{multline*}
 \frac{d^2}{2}\le\dep_{\Vv_\phi}\widetilde{(\phi^2)}_{\xi_g}
 =d\bigl(d-\deg_{\xi_g}(\phi)\bigr)
 +m_{\Vv_\phi}(\phi)\cdot\dep_{\Vv_\phi}\tilde{\phi}_{\xi_g}\\
 \le d\bigl(d-\deg_{\xi_g}(\phi)\bigr)+\frac{d-3}{2}\cdot m_{\Vv_\phi}(\phi).
\end{multline*}
Hence if $d<12$, then using the sharper estimate \eqref{eq:depfix} again, we further compute
\begin{multline*}
 m_{\Vv_\phi}(\phi)
 \ge\frac{d(-d+2\deg_{\xi_g}(\phi))}{d-3}
 =-d+2\deg_{\xi_g}(\phi)
 +3\frac{d-2\dep_{\Vv_\phi}\tilde{\phi}_{\xi_g}}{d-3}\\
 \ge-d+2\deg_{\xi_g}(\phi)+3\frac{d-(d-3)}{d-3}
 >-d+2\deg_{\xi_g}(\phi)+1.
\end{multline*}
\item Hence when $n=2$ and $d<12$, also using 
the injective morphism $\iota_e$ in \eqref{eq:coeff} and
Lemma \ref{th:ramif}, we have
\begin{align*}
 &\dim\bigl\{(\tilde{\phi}_{\xi_g},\dep_{\,\cdot}\tilde{\phi}_{\xi_g})
 \in C_{d,2}:\MHRL(\phi)=\{\xi_g\}\bigr\}\\
 \le&\dim\bigl\{\bigl(\Vv_\phi,(\tilde{\phi}_{\xi_g},\dep_{\,\cdot}\tilde{\phi}_{\xi_g})\bigr):
 (\tilde{\phi}_{\xi_g},\dep_{\,\cdot}\tilde{\phi}_{\xi_g})
 \in C_{d,2}\text{ and }\MHRL(\phi)=\{\xi_g\}\bigr\}\\
\le&\max\biggl\{d,\dim\bigcup_{1\le m\le e\le d: \atop m> -d+2e+1}
\bigl\{\bigl(\Vv,(\tilde{\phi}_{\xi_g},(d-e)1_{\Vv})\bigr)\in
\cR_{e,m}':\tilde{\phi}_{\xi_g}(\Vv)=\Vv\bigr\}\biggr\}\\
=&\max\biggl\{d,\dim\bigcup_{1\le m\le e\le d: \atop m\ge -d+2e+2}
\bigl\{(\Vv,\tilde{\phi}_{\xi_g})\in
\cR_{e,m}:\tilde{\phi}_{\xi_g}(\Vv)=\Vv\bigr\}\biggr\}\\
<&\max\Bigl\{d,\max_{1\le e\le d}\bigl(2e+3-(-d+2e+2)\bigr)\Bigr\}=d+1,
\end{align*}
which with \eqref{eq:fibers} and Lemma \ref{th:stab} yields
\begin{gather*}
 \dim\bigl(q_{d,k}(C_{d,2})\bigr)<(d+1)-3=d-2.  
\end{gather*}
\item 
Let us see the reverse estimate 
$\dim(q_{d,k}(C_{d,2}))\ge d-3$ when $d>7$.
We first claim that
$\dim\cR_{(d+3)/2,5}''=d$; for, the fiber over each 
$\Vv\in T_{\xi_g}\sP^1$
of the projection $\cR_{(d+3)/2,5}''\to T_{\xi_g}\sP^1$
is nothing but the fiber over $(\Vv,\Vv)\in\Delta_{T_{\xi_g}\sP^1}$
of the morphism $\cR_{(d+3)/2,5}\to T_{\xi_g}\sP^1\times T_{\xi_g}\sP^1$,
$(\Vv,\tilde{\phi}_{\xi_g})\mapsto(\Vv,\tilde{\phi}_{\xi_g}(\Vv))$,
in the proof of Lemma \ref{th:ramif},
so in particular is of dimension $2\cdot(d+3)/2+1-5=d-1$, 
and the claim holds. 
Next, for every 
$(\Vv,\tilde{\phi}_{\xi_g})\in\cR_{(d+3)/2,5}''$, pick
$\phi\in\Rat_{d,K}$ such that
$(\tilde{\phi}_{\xi_g},\dep_{\,\cdot}\tilde{\phi}_{\xi_g})
=\lambda_{(d+3)/2,5}(\Vv,\tilde{\phi}_{\xi_g})$,
which is not in $I(d)$ since $\deg_{\xi_g}(\phi)=(d+3)/2>0$.
Then we have not only $\MHRL(\phi)=\{\xi_g\}$
since by the slope formula \eqref{eq:slope}, 
$(d-1)\rd_{\Vv}\hypRes_\phi=-(d-3)/2+(d-1)/2>0$ 
and $(d-1)\rd_{\Vw}\hypRes_\phi\ge -0+(d-1)/2>0$
for every $\Vw\in(T_{\xi_g}\sP^1)\setminus\{\Vv\}$,
but also
\begin{align*}
 \dep_{\Vv}\widetilde{(\phi^2)}_{\xi_g}
 &=(\deg\phi)\cdot\dep_{\Vv}\tilde{\phi}_{\xi_g}
 +m_{\Vv}(\phi)\cdot
 \dep_{\tilde{\phi}_{\xi_g}(\Vv)}\tilde{\phi}_{\xi_g}\\
 &=d\cdot\frac{d-3}{2}+5\cdot\frac{d-3}{2}
 =\frac{(d-3)(d+5)}{2}\quad\text{and in turn}\\
(d^2-1)\rd_{\Vv}\hypRes_{\phi^2}
 &=-\dep_{\Vv}\widetilde{(\phi^2)}_{\xi_g}+\frac{d^2-1}{2}
 =-\frac{(d-3)(d+5)}{2}+\frac{d^2-1}{2}
 =7-d<0
\end{align*}
by the composition formula \eqref{eq:comp} applied to $\phi^2=\phi\circ\phi$ and the slope formula \eqref{eq:slope} applied to $\phi^2$.

Hence $\lambda_{(d+3)/2,5}(\cR_{(d+3)/2,5}'')\subset
\{(\tilde{\phi}_{\xi_g},\dep_{\,\cdot}\tilde{\phi}_{\xi_g})
\in C_{d,2}:\MHRL(\phi)=\{\xi_g\}\}$,
which with \eqref{eq:fibers} yields
\begin{gather*}
 \dim q_{d,k}(C_{d,2})
 \ge\dim\cR_{(d+3)/2,5}''-3=d-3,
\end{gather*}
and we are done in two of the five exceptional cases.
\end{inparaenum}

\item Suppose still $n=2$ but now $d\in\{5,7\}$.
\begin{inparaenum}[({ii}-1)]
\makeatletter
\renewcommand{\p@enumii}{\theenumi-}
\makeatother
\item\label{item:depbound} We have $\tilde{\phi}_{\xi_g}(\Vv_\phi)\ne\Vv_\phi$;
for, otherwise,
by Lemma \ref{th:depth},
the composition formula \eqref{eq:comp}
applied to $\phi^2=\phi\circ\phi$,
and the sharper estimate \eqref{eq:depfix}, we have
\begin{multline*}
 \frac{d^2}{2}\le\dep_{\Vv_\phi}\widetilde{(\phi^2)}_{\xi_g}
 =(\deg\phi)\cdot\dep_{\Vv_\phi}\tilde{\phi}_{\xi_g}
 +m_{\Vv_\phi}(\phi)\cdot\dep_{\tilde{\phi}_{\xi_g}(\Vv_\phi)}\tilde{\phi}_{\xi_g}\\
 \le d^2-d^2+d\cdot\dep_{\Vv_\phi}\tilde{\phi}_{\xi_g}
 +(d-\dep_{\Vv_\phi}\tilde{\phi}_{\xi_g})\cdot\dep_{\Vv_\phi}\tilde{\phi}_{\xi_g}
 =d^2-\bigl(d-\dep_{\Vv_\phi}\tilde{\phi}_{\xi_g}\bigr)^2\\
 \le d^2-\Bigl(d-\frac{d-3}{2}\Bigr)^2
 =\begin{cases}
 9&\text{if }d=5\\
 24&\text{if }d=7
 \end{cases}
 <\frac{d^2}{2},
\end{multline*}
which is impossible. Then using \eqref{eq:dep}, we further compute
\begin{multline}
 \Bigl\lceil\frac{d^2}{2}\Bigr\rceil
 \le\dep_{\Vv_\phi}\widetilde{(\phi^2)}_{\xi_g}
 =(\deg\phi)\cdot\dep_{\Vv_\phi}\tilde{\phi}_{\xi_g}
 +m_{\Vv_\phi}(\phi)\cdot
 \dep_{\tilde{\phi}_{\xi_g}(\Vv_\phi)}\tilde{\phi}_{\xi_g}\\
 \le d\cdot\dep_{\Vv_\phi}\tilde{\phi}_{\xi_g}
 +\frac{\bigl(m_{\Vv_\phi}(\phi)
 +\dep_{\tilde{\phi}_{\xi_g}(\Vv_\phi)}\tilde{\phi}_{\xi_g}\bigr)^2}{4}
 \le d\cdot\dep_{\Vv_\phi}\tilde{\phi}_{\xi_g}
 +\frac{\bigl(d-\dep_{\Vv_\phi}\tilde{\phi}_{\xi_g}\bigr)^2}{4}\\
 =\frac{\bigl(d+\dep_{\Vv_\phi}\tilde{\phi}_{\xi_g}\bigr)^2}{4}
\le \frac{\bigl(d+(d-1)/2)^2}{4}
 =\begin{cases}
 \dfrac{49}{4}&\text{if }d=5,\\
 25&\text{if }d=7,
 \end{cases}\label{eq:depbound}
\end{multline}
which is impossible when $d=5$.
Hence by Lemma \ref{th:stab}, we have $\dim(q_{5,k}(C_{5,2}))=0$,
and we are done in the exceptional case
\begin{gather*}
 (n,d)=(2,5). 
\end{gather*}
\item\label{item:cycle} 
Hence suppose in addition $d=7$. We have the injective morphism
\begin{gather}
 \mu_{2,7}:\cP_2\cap\cR_{2,2}\to\IRD_{7,k}(\xi_g),\quad
 (\Vv,\tilde{\phi}_{\xi_g})\mapsto
 \bigl(\tilde{\phi}_{\xi_g},3\cdot1_{\{\Vv\}}
 +2\cdot1_{\{\tilde{\phi}_{\xi_g}(\Vv)\}}\bigr)\label{eq:PR}
\end{gather}
from $\iota_2$ in \eqref{eq:coeff}. Let us see the inclusion
$\{(\tilde{\phi}_{\xi_g},\dep_{\,\cdot}\tilde{\phi}_{\xi_g})\in C_{d,n}:\MHRL(\phi)=\{\xi_g\}\}\subset\mu_{2,7}(\cP_2\cap\cR_{2,2})$;
indeed, in (\ref{item:depbound}),
all the equalities in \eqref{eq:depbound} are the cases,
 so in particular
 $(7+\dep_{\Vv_\phi}\tilde{\phi}_{\xi_g})^2/4=25$
 and $m_{\Vv_\phi}(\phi)=\dep_{\tilde{\phi}_{\xi_g}(\Vv_\phi)}\tilde{\phi}_{\xi_g}$, and then 
 $\dep_{\Vv_\phi}\tilde{\phi}_{\xi_g}=3$ and
 $m_{\Vv_\phi}(\phi)
 =\dep_{\tilde{\phi}_{\xi_g}(\Vv_\phi)}\tilde{\phi}_{\xi_g}=\sqrt{25-7\cdot3}=2$, and in turn $\deg_{\xi_g}(\phi)=2$ since
 $2=m_{\Vv_\phi}(\phi)\le\deg_{\xi_g}(\phi)\le d
 -\dep_{\tilde{\phi}_{\xi_g}(\Vv_\phi)}\tilde{\phi}_{\xi_g}-\dep_{\Vv_\phi}\tilde{\phi}_{\xi_g}=7-2-3=2$. Consequently
 \begin{gather*}
 \dep_{\,\cdot}\tilde{\phi}_{\xi_g}
 =3\cdot1_{\{\Vv_\phi\}}+2\cdot 1_{\{\tilde{\phi}_{\xi_g}(\Vv_\phi)\}}.
 \end{gather*}
 Moreover, noting that $\xi_g\notin\MHRL(\phi^2)$, by the slope formula \eqref{eq:slope} applied to $\phi^2$, 
 \begin{gather*}
 0>(7^2-1)\rd_{\Vw}\hypRes_{\phi^2}
 =-\dep_{\Vw}\widetilde{(\phi^2)}_{\xi_g}+\frac{1}{2}
 \begin{cases}
 7^2-1&\text{if }(\tilde{\phi}_{\xi_g})^2(\Vw)=\Vw,\\
 7^2+1&\text{otherwise},
 \end{cases}
 \end{gather*}
so that $\dep_{\Vw}\widetilde{(\phi^2)}_{\xi_g}\ge 25$
for some $\Vw\in T_{\xi_g}\sP^1$,
and then we have 
$\Vw=\Vv_\phi$, $\dep_{\Vv_\phi}\widetilde{(\phi^2)}_{\xi_g}=25$,
 and $(\tilde{\phi}_{\xi_g})^2(\Vv_\phi)=\Vv_\phi$
 also by Lemma \ref{th:depth} and \eqref{eq:depbound}. Hence
 $(\Vv_\phi,\tilde{\phi}_{\xi_g})$ is in $\cP_2\cap\cR_{2,2}$,
 and 
\begin{gather*}
\mu_{2,7}\bigl((\Vv_\phi,\tilde{\phi}_{\xi_g})\bigr)
=(\tilde{\phi}_{\xi_g},\dep_{\,\cdot}\tilde{\phi}_{\xi_g}).
\end{gather*}
We also claim the opposite inclusion 
$\{(\tilde{\phi}_{\xi_g},\dep_{\,\cdot}\tilde{\phi}_{\xi_g})\in C_{d,n}: \MHRL(\phi)=\{\xi_g\}\}\supset\mu_{2,7}(\cP_2\cap\cR_{2,2})$;
indeed, for every $(\Vv,\tilde{\phi}_{\xi_g})\in\cP_2\cap\cR_{2,2}$,
pick $\phi\in\Rat_{7,K}$ such that
$(\tilde{\phi}_{\xi_g},\dep_{\,\cdot}\tilde{\phi}_{\xi_g})
=\bigl(\tilde{\phi}_{\xi_g},3\cdot1_{\{\Vv\}}+2\cdot1_{\{\tilde{\phi}_{\xi_g}(\Vv)\}}\bigr)$,
which is not in $I(7)$ since $\deg_{\xi_g}(\phi)=2>0$.
Moreover, we have not only $\MHRL(\phi)=\{\xi_g\}$
since by the slope formula \eqref{eq:slope},
$(7-1)\rd_{\Vv}\hypRes_\phi=-3+(7+1)/2>0$
and $(7-1)\rd_{\Vu}\hypRes_\phi\ge-2+(7-1)/2>0$
for every $\Vu\in(T_{\xi_g}\sP^1)\setminus\{\Vv\}$,
but also $\xi_g\not\in\MHRL(\phi^2)$ since
\begin{gather*}
 \dep_{\Vv}\widetilde{(\phi^2)}_{\xi_g}
 =(\deg\phi)\cdot\dep_{\Vv}\tilde{\phi}_{\xi_g}
 +m_{\Vv}(\phi)\cdot
 \dep_{\tilde{\phi}_{\xi_g}(\Vv)}\tilde{\phi}_{\xi_g}
 =7\cdot3+2\cdot2=25\quad\text{and in turn}\\
 (7^2-1)\rd_{\Vv}\hypRes_{\phi^2}
 =-\dep_{\Vv}\widetilde{(\phi^2)}_{\xi_g}+\frac{7^2-1}{2}
 =-25+24=-1<0
\end{gather*}
by the composition formula \eqref{eq:comp} applied to $\phi^2=\phi\circ\phi$ and the slope formula \eqref{eq:slope} applied to $\phi^2$. Hence the claim holds.
Now also by \eqref{eq:fibers} and Lemmas \ref{th:stab} and \ref{th:dimcycle},
\begin{gather*}
 \dim q_{7,k}(C_{7,2})=\dim(\cP_2\cap\cR_{2,2})-3=4-3=1,
\end{gather*}
and we are done in the exceptional case
\begin{gather*}
 (n,d)=(2,7). 
\end{gather*}
\end{inparaenum}

\item Suppose finally that $n=3$ and $d=5$.
\begin{inparaenum}[({iii}-1)]
\makeatletter
\renewcommand{\p@enumii}{\theenumi-}
\makeatother
\item By Lemma \ref{th:depth} and the composition formula \eqref{eq:comp}
applied to $\phi^3=\phi^2\circ\phi$ and then to $\phi^2=\phi\circ\phi$, we have
\begin{multline}
 \frac{d^3}{2}\le\dep_{\Vv_\phi}\widetilde{(\phi^3)}_{\xi_g}
 =\bigl(\deg(\phi^2)\bigr)\cdot\dep_{\Vv_\phi}\tilde{\phi}_{\xi_g}
 +m_{\Vv_\phi}(\phi)\cdot
 \dep_{\tilde{\phi}_{\xi_g}(\Vv_\phi)}\widetilde{(\phi^2)}_{\xi_g}\\
 =25\cdot\dep_{\Vv_\phi}\tilde{\phi}_{\xi_g}
 +m_{\Vv_\phi}(\phi)\cdot\Bigl(
 (\deg\phi)\cdot\dep_{\tilde{\phi}_{\xi_g}(\Vv_\phi)}\tilde{\phi}_{\xi_g}
 +m_{\tilde{\phi}_{\xi_g}(\Vv_\phi)}(\phi)
 \cdot\dep_{(\tilde{\phi}_{\xi_g})^2(\Vv_\phi)}\tilde{\phi}_{\xi_g}\Bigr).
 \label{eq:cub}
\end{multline}
Let us see 
\begin{gather*}
 \dep_{\Vv_\phi}\tilde{\phi}_{\xi_g}=2; 
\end{gather*}
for, otherwise, we have 
$\dep_{\,\cdot}\tilde{\phi}_{\xi_g}=1_{\{\Vv_\phi\}}$ since
$0\le(5-m_{\Vv_\phi}(\phi))/2<\dep_{\Vv_\phi}\tilde{\phi}_{\xi_g}<(5+1)/2=3$ by \eqref{eq:dep} so $\dep_{\Vv_\phi}\tilde{\phi}_{\xi_g}=1$
and $3<m_{\Vv_\phi}(\phi)\le\deg_{\xi_g}(\phi)
\le 5-\dep_{\Vv_\phi}\tilde{\phi}_{\xi_g}=4$, and in turn
also by \eqref{eq:cub},
$125/2=d^3/2\le
 5^2\cdot 1+4\cdot\bigl(5\cdot 1+\deg_{\xi_g}(\phi)\cdot 1\bigr)
=25+4(5+4\cdot 1)=61$,
which is impossible. Then we moreover have not only
\begin{gather*}
 \tilde{\phi}_{\xi_g}(\Vv_\phi)\ne\Vv_\phi 
\end{gather*}
since otherwise,
by the sharper estimate \eqref{eq:depfix}, 
$2=\dep_{\Vv_\phi}\tilde{\phi}_{\xi_g}<(5-1)/2=2$,
which is impossible, but also $\deg_{\xi_g}(\phi)\in\{2,3\}$ 
since using \eqref{eq:dep} again, 
\begin{gather}
 1=\frac{\dep_{\Vv_\phi}\tilde{\phi}_{\xi_g}}{2}
\le\frac{d-\deg_{\xi_g}(\phi)}{2}
=\frac{5-\deg_{\xi_g}(\phi)}{2}
\le\frac{5-m_{\Vv_\phi}(\phi)}{2}<\dep_{\Vv_\phi}\tilde{\phi}_{\xi_g}=2.
\label{eq:deg}
\end{gather}
\item\label{item:cub} When $\deg_{\xi_g}(\phi)=3>0$, we have
$\dep_{\,\cdot}\tilde{\phi}_{\xi_g}=2\cdot 1_{\{\Vv_\phi\}}$.
Then also by \eqref{eq:cub}, we compute
$125/2\le
 25\cdot 2+m_{\Vv_\phi}(\phi)
(5\cdot 0+m_{\tilde{\phi}_{\xi_g}(\Vv_\phi)}(\phi)
\cdot\dep_{(\tilde{\phi}_{\xi_g})^2(\Vv_\phi)}\tilde{\phi}_{\xi_g})$, which yields $(\tilde{\phi}_{\xi_g})^2(\Vv_\phi)=\Vv_\phi$
and $m_{\Vv_\phi}(\phi)=m_{\tilde{\phi}_{\xi_g}(\Vv_\phi)}(\phi)=3$. 
Then there is $M\in\PGL(2,K^\circ)$ such that
$\tilde{M}_{\xi_g}(\Vv_\phi)=\widevec{\xi_g\infty}$ and
$\tilde{M}_{\xi_g}\bigl(\tilde{\phi}_{\xi_g}(\Vv_\phi)\bigr)=\widevec{\xi_g0}$ and that
$(\widetilde{(M\circ\phi\circ M^{-1})}_{\xi_g},
\dep_{\,\cdot}\widetilde{(M\circ\phi\circ M^{-1})}_{\xi_g})
=(\widetilde{(z^{-3})}_{\xi_g},2\cdot1_{\{\widevec{\xi_g\infty}\}})$
in $\IRD_{5,k}(\xi_g)$, so that
\begin{gather*}
 q_{5,k}\bigl((\tilde{\phi}_{\xi_g},\dep_{\,\cdot}\tilde{\phi}_{\xi_g})\bigr)
=q_{5,k}\Bigl(\bigl(\widetilde{(z^{-3})}_{\xi_g},2\cdot 1_{\{\widevec{\xi_g\infty}\}}\bigr)\Bigr).
\end{gather*}
Alternatively when $\deg_{\xi_g}(\phi)=2$, 
let us first see the inclusion
\begin{gather*}
 \bigl\{(\tilde{\phi}_{\xi_g},\dep_{\,\cdot}\tilde{\phi}_{\xi_g})
\in C_{5,3}:\MHRL(\phi)=\{\xi_g\}\text{ and }\tilde{\phi}_{\xi_g}\in\Rat_{2,k}\bigr\}
\subset\mu_{2,5}(\cP_2\cap\cR_{2,2}), 
\end{gather*}
using the injective morphism 
\begin{gather}
 \mu_{2,5}:\cP_2\cap\cR_{2,2}\to\IRD_{5,k}(\xi_g),\quad
 (\Vv,\tilde{\phi}_{\xi_g})\mapsto
 \bigl(\tilde{\phi}_{\xi_g},2\cdot1_{\{\Vv\}}
 +1_{\{\tilde{\phi}_{\xi_g}(\Vv)\}}\bigr)\label{eq:PC}
\end{gather}
from $\iota_2$ in \eqref{eq:coeff}; indeed, 
by \eqref{eq:deg}, we have 
\begin{gather*}
 m_{\Vv_\phi}(\phi)=2 
\end{gather*}
and, for some $\Vw\in(T_{\xi_g}\sP^1)\setminus\{\Vv_\phi\}$,
$\dep_{\,\cdot}\tilde{\phi}_{\xi_g}
=2\cdot1_{\{\Vv_\phi\}}+1_{\{\Vw\}}$. 
Then we also have
\begin{gather*}
 (\tilde{\phi}_{\xi_g})^2(\Vv_\phi)=\Vv_\phi
\quad\text{and}\quad \tilde{\phi}_{\xi_g}(\Vv_\phi)=\Vw; 
\end{gather*}
indeed, unless $\tilde{\phi}_{\xi_g}(\Vv_\phi)=\Vw$,
we have $\dep_{\tilde{\phi}_{\xi_g}(\Vv_\phi)}\tilde{\phi}_{\xi_g}=0$,
so that by \eqref{eq:cub}, 
$125/2\le 25\cdot 2+2(5\cdot 0
+\deg_{\xi_g}(\phi)\cdot\max\{2,1\})
\le 50+2\cdot 2\cdot 2=58$,
which is impossible, and in turn unless $\tilde{\phi}_{\xi_g}(\Vw)\neq\Vw$, we have $\tilde{\phi}_{\xi_g}(\Vv_\phi)=\tilde{\phi}_{\xi_g}(\Vw)=\Vw$ so that by \eqref{eq:directional}, $3=2+1\le m_{\Vv_\phi}(\phi)+m_{\Vw}(\phi)
\le \sum_{\Vu\in T_{\xi_g}\sP^1:\tilde{\phi}_{\xi_g}(\Vu)=\Vw}m_{\Vu}(\phi) =\deg_{\xi_g}(\phi)=2$, which is impossible, and then
unless $\tilde{\phi}_{\xi_g}(\Vw)=\Vv_\phi$, we have
$\dep_{\tilde{\phi}_{\xi_g}(\Vw)}\tilde{\phi}_{\xi_g}
=0$, so that by \eqref{eq:cub} again, 
$125/2\le 25\cdot 2+2\bigl(5\cdot\dep_{\Vw}\tilde{\phi}_{\xi_g}
+m_{\tilde{\phi}_{\xi_g}(\Vv_\phi)}(\phi)\cdot\dep_{\tilde{\phi}_{\xi_g}(\Vw)}\tilde{\phi}_{\xi_g}\bigr)
=50+2(5\cdot 1+0)=60$, which is impossible.
Hence $(\Vv_\phi,\tilde{\phi}_{\xi_g})$ is in $\cP_2\cap\cR_{2,2}$,
and 
\begin{gather*}
\mu_{2,5}\bigl((\Vv_\phi,\tilde{\phi}_{\xi_g})\bigr)
=(\tilde{\phi}_{\xi_g},\dep_{\,\cdot}\tilde{\phi}_{\xi_g}).
\end{gather*}
We also claim the opposite inclusion
\begin{gather*}
 \bigl\{(\tilde{\phi}_{\xi_g},\dep_{\,\cdot}\tilde{\phi}_{\xi_g})
\in C_{5,3}:\MHRL(\phi)=\{\xi_g\}\text{ and }\tilde{\phi}_{\xi_g}\in\Rat_{2,k}\bigr\}
\supset\mu_{2,5}(\cP_2\cap\cR_{2,2});
\end{gather*}
indeed, for every $(\Vv,\tilde{\phi}_{\xi_g})\in\cP_2\cap\cR_{2,2}$,
pick $\phi\in\Rat_{5,K}$ such that
$(\tilde{\phi}_{\xi_g},\dep_{\,\cdot}\tilde{\phi}_{\xi_g})
=\bigl(\tilde{\phi}_{\xi_g},2\cdot1_{\{\Vv\}}+1_{\{\tilde{\phi}_{\xi_g}(\Vv)\}}\bigr)$,
which is not in $I(5)$ since $\deg_{\xi_g}(\phi)=2>0$.
Moreover, we have not only $\MHRL(\phi)=\{\xi_g\}$
since by the slope formula \eqref{eq:slope},
$(5-1)\rd_{\Vv}\hypRes_\phi=-2+(5+1)/2>0$
and $(5-1)\rd_{\Vu}\hypRes_\phi\ge-1+(5-1)/2>0$
for every $\Vu\in(T_{\xi_g}\sP^1)\setminus\{\Vv\}$,
but also $\xi_g\not\in\MHRL(\phi^3)$ since, using
$(\tilde{\phi}_{\xi_g})^2(\Vv)=\Vv$,
\begin{gather*}
 \dep_{\Vv}\widetilde{(\phi^3)}_{\xi_g}
 =25\cdot 2+2\bigl(5\cdot\dep_{\Vw}\tilde{\phi}_{\xi_g}+m_{\tilde{\phi}_{\xi_g}(\Vv)}(\phi)\cdot
 \dep_{\Vv}\tilde{\phi}_{\xi_g}\bigr)
 \ge 50+2(5\cdot 1+1\cdot 2)=64\\
\text{and in turn}\quad
 (5^3-1)\rd_{\Vv}\hypRes_{\phi^3}
 \le-\dep_{\Vv}\widetilde{(\phi^3)}_{\xi_g}+\frac{5^3+1}{2}
 \le-64+63=-1<0
\end{gather*}
using \eqref{eq:cub} and by the slope formula \eqref{eq:slope}
applied to $\phi^3$.
Hence the claim holds.
Now also by \eqref{eq:fibers} and Lemmas \ref{th:stab} and \ref{th:dimcycle},
we have
$\dim q_{5,k}(C_{5,3})=\dim(\cP_2\cap\cR_{2,2})-3=4-3=1$,
and we are done in the exceptional case of
$(n,d)=(3,5)$.
\end{inparaenum}
\end{inparaenum}

Now the proof of Theorem \ref{th:main}(\ref{item:moduli}) is complete.

\section{Proofs of Theorem \ref{th:cubic} and Corollary \ref{th:indet}}\label{sec:cubicindet}

\subsection{Proof of Theorem \ref{th:cubic}}

Let $d=3$ and $K$ be as in Section \ref{sec:intro}.

For every $n\ge 2$, by Lemma \ref{th:stab},
it remains to show that for every $\phi\in\Rat_{3,K}$ such that
$(\tilde{\phi}_{\xi_g},\dep_{\,\cdot}\tilde{\phi}_{\xi_g})\in C_{3,n}$
but $\MHRL(\phi)=\{\xi_g\}$, there is $M\in\PGL(2,K^\circ)$ such that
\begin{gather*}
 \Bigl(
 \widetilde{\bigl(M\circ\phi\circ M^{-1}\bigr)}_{\xi_g},
 \dep_{\,\cdot}\widetilde{\bigl(M\circ\phi\circ M^{-1}\bigr)}_{\xi_g}
 \Bigr)
 =\Bigl(
 \widetilde{(z^{-2})}_{\xi_g},
 1_{\{\widevec{\xi_g\infty}\}}
 \Bigr)\quad\text{in }\IRD_{3,k}(\xi_g).
\end{gather*}
Indeed, for any such $\phi\in\Rat_{3,K}$,
\begin{gather*}
 \phi(\xi_g)=\xi_g
\end{gather*}
since otherwise, by Lemma \ref{th:depth} and
the slope formula \eqref{eq:slope}, we must have
$3/2=d/2<\dep_{\Vv_\phi}\tilde{\phi}_{\xi_g}<(d+1)/2=2$,
which is impossible. 
Then by Lemma \ref{th:depth} and
the slope formula \eqref{eq:slope} again, we have
$0<2\dep_{\Vv_\phi}\tilde{\phi}_{\xi_g}<d+1=4$, that is,
\begin{gather*}
 \dep_{\Vv_\phi}\tilde{\phi}_{\xi_g}=1\le\frac{d}{2}, 
\end{gather*}
and then we also have not only 
\begin{gather*}
 \tilde{\phi}_{\xi_g}(\Vv_\phi)\neq\Vv_\phi 
\end{gather*}
since otherwise the impossible $0<-1+(d-1)/2=0$ occurs,
but also
$1=3-2\cdot 1=d-2\dep_{\Vv_\phi}\tilde{\phi}_{\xi_g}
<m_{\Vv_\phi}(\phi)\le\deg_{\xi_g}(\phi)\le d-\dep_{\Vv_\phi}\tilde{\phi}_{\xi_g}=3-1=2$ so 
\begin{gather*}
 m_{\Vv_\phi}(\phi)=\deg_{\xi_g}(\phi)=2.
\end{gather*}
Then we also compute
\begin{gather*}
 \sum_{\Vw\in(T_{\xi_g}\sP^1)\setminus\{\Vv_\phi\}}\dep_{\Vw}\tilde{\phi}_{\xi_g}=d-\deg_{\xi_g}(\phi)-\dep_{\Vv_\phi}\tilde{\phi}_{\xi_g}=3-2-1=0.
\end{gather*}
Hence by Lemma \ref{th:depth}, $\Vv_\phi$ is periodic 
under $\tilde{\phi}_{\xi_g}$ having the exact period $p$ in
$\{2,\ldots,n\}$, and also using the composition formula 
\eqref{eq:comp} repeatedly, we have
\begin{multline*}
\frac{3^n}{2}
=\frac{\deg(\phi^n)}{2}\le \dep_{\Vv_\phi}\widetilde{(\phi^n)}_{\xi_g}
=\deg(\phi^{n-1})\cdot\dep_{\Vv_\phi}\tilde{\phi}_{\xi_g}+m_{\Vv_\phi}(\phi)\cdot\dep_{\tilde{\phi}_{\xi_g}(\Vv_\phi)}\widetilde{(\phi^{n-1})}_{\xi_g}\\
 =\sum_{j=0}^{n-1}m_{\Vv_\phi}(\phi^{n-1-j})\cdot 3^j\cdot 
 \dep_{(\tilde{\phi}_{\xi_g})^{n-1-j}(\Vv_\phi)}\tilde{\phi}_{\xi_g}
 =\sum_{\ell=0}^{\lfloor(n-1)/p\rfloor}
 3^{n-1-\ell p}\cdot (m_{\Vv_{\phi}}(\phi^p))^\ell\cdot\dep_{\Vv_\phi}\tilde{\phi}_{\xi_g},
\end{multline*}
so that $3/2\le\sum_{\ell=0}^{\lfloor(n-1)/p\rfloor}(m_{\Vv_{\phi}}(\phi^p)/3^p)^{\ell}<1/(1-(m_{\Vv_{\phi}}(\phi^p)/3^p))
\le 1/(1-(\deg_{\xi_g}(\phi)/3)^p)$. 
Hence we have
\begin{gather*}
 p=2\quad\text{and}\quad m_{\tilde{\phi}_{\xi_g}\Vv_\phi}(\phi)=2
\end{gather*}
by $(2/3)^p>1/3$ and in turn by $2\ge m_{\tilde{\phi}_{\xi_g}\Vv_\phi}(\phi)>3/2$ also using $m_{\Vv_\phi}(\phi)=2$. 
Now we are done.

Conversely, pick $\phi\in\Rat_{3,K}$ such that $(\tilde{\phi}_{\xi_g},\dep_{\,\cdot}\tilde{\phi}_{\xi_g})=(\widetilde{(z^{-2})}_{\xi_g},1_{\{\widevec{\xi_g\infty}\}})$. Then $\MHRL(\phi)=\{\xi_g\}$
since by the slope formula \eqref{eq:slope}, 
\begin{gather*}
 \rd_{\widevec{\xi_g\infty}}\hypRes_\phi
 =\frac{1}{3-1}\Bigl(-1+\frac{3+1}{2}\Bigr)=\frac{1}{2}>0
 \quad\text{and}\\
 \rd_{\Vw}\hypRes_\phi\ge\frac{1}{d-1}\Bigl(-0+\frac{d-1}{2}\Bigr)=\frac{1}{2}>0
 \quad\text{for every }\Vw\in(T_{\xi_g}\sP^1)\setminus\{\widevec{\xi_g\infty}\}.
\end{gather*}
Pick also any $n\ge 2$. Then
$\widetilde{(\phi^n)}_{\xi_g}=(\phi^n)_{*,\xi_g}
=(\phi_{*,\xi_g})^n=(\tilde{\phi}_{\xi_g})^n
=(\widetilde{(z^{-2})}_{\xi_g})^n$ and
$m_{\widevec{\xi_g\infty}}(\phi)=m_{\tilde{\phi}_{\xi_g}\widevec{\xi_g\infty}}(\phi)=2$
so $m_{\widevec{\xi_g\infty}}(\phi^2)=m_{\tilde{\phi}_{\xi_g}\widevec{\xi_g\infty}}(\phi)\cdot m_{\widevec{\xi_g\infty}}(\phi)=2^2$, and then 
by the composition formula \eqref{eq:comp},
we have $\dep_{\widevec{\xi_g\infty}}\widetilde{(\phi^2)}_{\xi_g}
=d\cdot\dep_{\widevec{\xi_g\infty}}\tilde{\phi}_{\xi_g}+2\cdot\dep_{\tilde{\phi}_{\xi_g}\widevec{\xi_g\infty}}\tilde{\phi}_{\xi_g}
=3\cdot 1+0=3$ and in turn recursively
\begin{multline*}
 \dep_{\widevec{\xi_g\infty}}\widetilde{(\phi^n)}_{\xi_g}
 =\deg(\phi^{n-2})\cdot\dep_{\widevec{\xi_g\infty}}\widetilde{(\phi^2)}_{\xi_g}
 +m_{\widevec{\xi_g\infty}}(\phi^2)\cdot\dep_{\widevec{\xi_g\infty}}\widetilde{(\phi^{n-2})}_{\xi_g}\\
 =3^{n-2}\cdot3+2^2\cdot\dep_{\widevec{\xi_g\infty}}\widetilde{(\phi^{n-2})}_{\xi_g} =\sum_{j=0}^{\lfloor(n-1)/2\rfloor}3^{n-1-2j}\cdot2^{2j}
 =3^{n-1}\sum_{j=0}^{\lfloor(n-1)/2\rfloor}\Bigl(\frac{2}{3}\Bigr)^{2j}.
\end{multline*}
When $n\in\{2,3,4\}$, we have $\xi_g\in\MHRL(\phi^n)$ and then 
$(\tilde{\phi}_{\xi_g},\dep_{\,\cdot}\tilde{\phi}_{\xi_g})\not\in C_{3,n}$; for, otherwise,
also by Lemma \ref{th:depth}, we have $\Vv_\phi=\widevec{\xi_g\infty}$ 
and moreover we must have
\begin{gather*}
 \frac{3^n}{2}\le\dep_{\widevec{\xi_g\infty}}\widetilde{(\phi^n)}_{\xi_g}
 \le 3^{n-1}\bigl(1+\Bigl(\frac{2}{3}\Bigr)^2\bigr)
 =13\cdot3^{n-3}<\frac{3^n}{2},
\end{gather*}
which is a contradiction. When $n\ge 5$, we have $\xi_g\not\in\MHRL(\phi^n)$, and then $(\widetilde{(z^{-2})}_{\xi_g},1_{\{\widevec{\xi_g\infty}\}})\in C_{3,n}$ and $q_{3,k}(\widetilde{(z^{-2})}_{\xi_g},1_{\{\widevec{\xi_g\infty}\}})\neq q_{3,k}(\widevec{\xi_g\infty},1_{\{\widevec{\xi_g\infty}\}}+2\cdot 1_{\{\widevec{\xi_g0}\}})$ (see the first component of \eqref{eq:stable}),
since also by the slope formula \eqref{eq:slope} applied to $\phi^n$, we have
\begin{multline*}
 (3^n-1)\rd_{\widevec{\xi_g\infty}}\hypRes_{\phi^n}
 \le -\Bigl(3^{n-1}\bigl(1+\Bigl(\frac{2}{3}\Bigr)^2+\Bigl(\frac{2}{3}\Bigr)^4+0\bigr)\Bigr)+\frac{3^n+1}{2}\\
 =-133\cdot3^{n-5}+\frac{3^n+1}{2}
 =\frac{1-23\cdot3^{n-5}}{2}<0.
\end{multline*}

Now the proof of Theorem \ref{th:cubic} is complete. \qed

\subsection{Proof of Corollary \ref{th:indet}}\label{sec:indet}

\begin{lemma}\label{th:indetdim}
For every $d\ge 2$, we have $\dim I(d)=d$.
When $d>2$, we have
\begin{gather*}
 \dim q_{d,k}\bigl(\{(\tilde{\phi}_{\xi_g},\dep_{\,\cdot}\tilde{\phi}_{\xi_g})\in I(d):\xi_g\in\MHRL(\phi)\}\bigr)\le d-3,
\end{gather*}
and when $d=2$,
$\{(\tilde{\phi}_{\xi_g},\dep_{\,\cdot}\tilde{\phi}_{\xi_g})\in I(2):\xi_g\in\MHRL(\phi)\}=\emptyset$.
\end{lemma}

\begin{proof}
The morphism $I(d)\to T_{\xi_g}\sP^1$,
$(\Vw,\dep_{\,\cdot}\tilde{\phi}_{\xi_g})\mapsto\Vw$,
has fibers isomorphic to $\bP_k^{d-1}$,
so $\dim I(d)=1+(d-1)=d$.

Pick any $\phi\in\Rat_{d,K}$ such that $(\tilde{\phi}_{\xi_g},\dep_{\,\cdot}\tilde{\phi}_{\xi_g})\in I(d)$
and that $\xi_g\in\MHRL(\phi)$.
Let us see
\begin{gather*}
 \#\bigl\{\Vv\in T_{\xi_g}\sP^1:\dep_{\Vv}\tilde{\phi}_{\xi_g}>0\bigr\}\ge 2;
\end{gather*}
for, if there is a unique $\Vv\in T_{\xi_g}\sP^1$ such that 
$\dep_{\Vv}\tilde{\phi}_{\xi_g}>0$, 
then $\tilde{\phi}_{\xi_g}\equiv\Vv$ and $\dep_{\Vv}\tilde{\phi}_{\xi_g}=d$, and in turn by $\xi_g\in\MHRL(\phi)$ and the slope formula \eqref{eq:slope}, we must have
$0\le(d-1)\rd_{\Vv}\hypRes_\phi=(-d+(d-1)/2)<0$,
which is a contradiction.

For any $M\in\PGL(2,K^\circ)$ fixing $(\tilde{\phi}_{\xi_g},\dep_{\,\cdot}\tilde{\phi}_{\xi_g})$ under the action \eqref{eq:stable} of $\PGL(2,K^\circ)$ on $\IRD_{d,k}(\xi_g)$, $\tilde{M}_{\xi_g}$ acts on
the finite set $\{\Vv\in T_{\xi_g}\sP^1:\dep_{\Vv}\tilde{\phi}_{\xi_g}>0\}$. Hence if $\#\bigl\{\Vv\in T_{\xi_g}\sP^1:\dep_{\Vv}\tilde{\phi}_{\xi_g}>0\bigr\}>2$, which is never the case when $d=2$,
then recalling Subsection \ref{sec:strat},
the stabilizer subgroup in $\PGL(2,K^\circ)$ or strictly speaking in $\PGL(2,K^\circ)/\bigl\{M\in\PGL(2,K^\circ):\tilde{M}_{\xi_g}=\mathrm{id}\bigr\}$ of $(\tilde{\phi}_{\xi_g},\dep_{\,\cdot}\tilde{\phi}_{\xi_g})$ is finite.
Alternatively if $\#\bigl\{\Vv\in T_{\xi_g}\sP^1:\dep_{\Vv}\tilde{\phi}_{\xi_g}>0\bigr\}=2$, then recalling $(\tilde{\phi}_{\xi_g},\dep_{\,\cdot}\tilde{\phi}_{\xi_g})\in I(d)$, we write
$\bigl\{\Vv\in T_{\xi_g}\sP^1:\dep_{\Vv}\tilde{\phi}_{\xi_g}>0\bigr\}=\{\widevec{\xi_g\phi(\xi_g)},\Vw\}$ so that
$\dep_{\widevec{\xi_g\phi(\xi_g)}}\tilde{\phi}_{\xi_g}+\dep_{\Vw}\tilde{\phi}_{\xi_g}=d$, and also recalling $\xi_g\in\MHRL(\phi)$,
by the slope formula \eqref{eq:slope}, we also have
\begin{gather*}
 \dep_{\widevec{\xi_g\phi(\xi_g)}}\tilde{\phi}_{\xi_g}\le\frac{d-1}{2}
\quad\text{and}\quad
 \dep_{\Vw}\tilde{\phi}_{\xi_g}
 \le\frac{d+1}{2},
\end{gather*}
so that the equalities are the case, that $2\nmid d$ so $d\ge 3$,
and that for some $M\in\PGL(2,K^\circ)$,
\begin{gather*}
 \Bigl(\widetilde{(M\circ\phi\circ M^{-1})}_{\xi_g},
 \dep_{\,\cdot}\widetilde{(M\circ\phi\circ M^{-1})}_{\xi_g}\Bigr)
 =\Bigl(\widevec{\xi_g\infty},
 \frac{d-1}{2}1_{\{\widevec{\xi_g\infty}\}}
 +\frac{d+1}{2}1_{\{\widevec{\xi_g0}\}}\Bigr)
\quad\text{in }\IRD_{d,k}(\xi_g).
\end{gather*}
Consequently, when $d>2$, we have
\begin{gather*}
 \dim q_{d,k}\bigl(\{(\tilde{\phi}_{\xi_g},\dep_{\,\cdot}\tilde{\phi}_{\xi_g})\in I(d):\xi_g\in\MHRL(\phi)\}\bigr)\le\max\{\dim I(d)-\dim\PGL(2,k),0\}\le d-3,
\end{gather*}
and when $d=2$,
$\{(\tilde{\phi}_{\xi_g},\dep_{\,\cdot}\tilde{\phi}_{\xi_g})\in I(2):\xi_g\in\MHRL(\phi)\}=\emptyset$.
\end{proof}

By Theorem \ref{th:moduli}, we have
$\dim I(\Phi_{d,n,k})=d-2$ if and only if
$\dim q_{d,k}(C_{d,n})=d-2$.
Now the proof of Corollary \ref{th:indet} is complete
by Theorems
\ref{th:main}(\ref{item:moduli}) and \ref{th:cubic}. \qed

\appendix
\refstepcounter{section}
\section*{Appendix: Indeterminacy of iterate maps}\label{app:coeff}\label{sec:modindet}

\begin{applemma}\label{th:graded}
Let $k$ be a field, and 
let $\theta:S\to R$ be a $k$-algebra homomorphism between
finitely generated nonnegatively graded integral $k$-algebras $R,S$
such that
\begin{gather*}
 R_0=S_0=k
\end{gather*}
and that
for some integer $r\ge 1$, we have $\theta(S_m)\subset R_{rm}$ 
for every $m\ge 0$.
Set $S_+=\bigoplus_{m>0}S_m$ and assume that 
$\theta(S_+)\ne 0$ and that $\operatorname{codim}_{\Proj R}V_+\bigl(R\theta(S_+)\bigr)\ge 2$.
Then the indeterminacy locus of
the rational map 
$\Proj R\dashrightarrow\Proj S$
induced by $\theta$ coincides with
$V_+\bigl(R\theta(S_+)\bigr)$.
\end{applemma}

\begin{proof}
Fix $\ell\gg 1$ such that the $\ell$th Veronese subrings
of $R$ and $S$ are generated by $R_\ell$ and $S_\ell$,
respectively, and fix a basis $s_0,\ldots,s_N$ of $S_\ell$,
so that we have the closed embedding
$[s_0:\cdots:s_N]:\Proj S\hookrightarrow\bP_k^N$.
By \cite[II, Theorem 7.1]{Hartshorne}, we have the morphism
$[\theta(s_0):\cdots:\theta(s_N)]:
\Proj R\setminus V_+\bigl(R\theta(S_+)\bigr)
\to\Proj S\hookrightarrow\bP_k^N$.
If this map extends regularly at some point of
$V_+\bigl(R\theta(S_+)\bigr)$,
then, after a local trivialization of
$\mathcal{O}_{\Proj R}(r\ell)$,
there is some $i\in\{0,\ldots,N\}$ such that
$\theta(s_j)/\theta(s_i)$ is regular at that point
for every $j\in\{0,\ldots,N\}$,
so the ideal generated by $\theta(s_0),\ldots,\theta(s_N)$
in the local ring at that point is the non-zero principal ideal
$(\theta(s_i))$. This contradicts the codimension hypothesis by
Krull's Hauptidealsatz \cite[I, Theorem 1.11A]{Hartshorne}.
\end{proof}

Let $d>1$ and $K$ be as in Section \ref{sec:intro}.
We show a statement slightly better than the following.

\begin{apptheorem}\label{th:coeff}
For every $n\ge 2$, the indeterminacy loci 
$I(\Phi^{\mathrm{IRD}}_{d,n,k})$ and $I(\Phi_{d,n,k})$
of $\Phi^{\mathrm{IRD}}_{d,n,k}$ and $\Phi_{d,n,k}$
respectively coincide with 
$I(d)$ and
$q_{d,k}\bigl(\bigl\{(\tilde{\phi}_{\xi_g},\dep_{\,\cdot}\tilde{\phi}_{\xi_g})\in I(d):\xi_g\in\MHRL(\phi)\bigr\}\cup C_{d,n}\bigr)$.
\end{apptheorem}

\begin{proof}
Pick any $n\ge 2$. Regarding $V_{d,k}$ as $k[X,Y]_d^{\oplus 2}$,
we have the $n$th iterate map
\begin{gather*}
 \Phi_n:V_{d,k}\to V_{d^n,k},\quad F\mapsto F^{\circ n},
\end{gather*}
which descends to
the rational map
$\Phi^{\mathrm{IRD}}_{d,n,k}:
\IRD_{d,k}(\xi_g)\dashrightarrow\IRD_{d^n,k}(\xi_g)$
using \eqref{eq:ird}, and for every $m\ge 0$, we have
$\Phi_n^*\bigl((\Sym(V_{d^n,k}^{\vee}))_m\bigr)
 \subset(\Sym(V_{d,k}^{\vee}))_{(1+d+\cdots+d^{n-1})m}$.
For every $\phi\in\Rat_{d,K}$, every
$F\in V_{d,k}\setminus\{0\}$ such that $[F]=\hat{\phi}$,
and every $j\ge 0$ such that $F^{\circ j}\ne 0$,
recalling \eqref{eq:ird} and \eqref{eq:coeff}, we have
$F^{\circ(j+1)}=0$ in $V_{d^{j+1},k}$ if and only if
$\phi^j(\xi_g)\ne\xi_g$ and
$\dep_{\widevec{\xi_g\phi^j(\xi_g)}}\tilde{\phi}_{\xi_g}>0$.
We claim that
\begin{gather}
 V_+\bigl(\Phi_n^*(V_{d^n,k}^{\vee})\bigr)=I(d);\label{eq:indet}
\end{gather}
indeed, if $\phi(\xi_g)=\xi_g$, then 
we have $F^{\circ n}\in V_{d^n,k}\setminus\{0\}$
using the above equivalence inductively, and 
if $\phi(\xi_g)\ne\xi_g$ and
$\dep_{\widevec{\xi_g\phi(\xi_g)}}\tilde{\phi}_{\xi_g}=0$,
then as seen in the proof of Lemma \ref{th:depth},
we have $\phi^j(\xi_g)\in U(\widevec{\xi_g\phi(\xi_g)})$ 
for every $j\ge 1$, and in turn still $F^{\circ n}\in V_{d^n,k}\setminus\{0\}$ using the above equivalence again inductively. 
In the remaining case of $(\tilde{\phi}_{\xi_g},\dep_{\,\cdot}\tilde{\phi}_{\xi_g})\in I(d)$, 
by the above equivalence for $j=1$,
$F^{\circ 2}=0$ in $V_{d^2,k}$, and in turn  
$F^{\circ n}=0$ in $V_{d^n,k}$.
Once \eqref{eq:indet} is at our disposal, we have 
$\operatorname{codim}_{\IRD_{d,k}(\xi_g)}V_+(\Phi_n^*(V_{d^n,k}^{\vee}))=(2d+1)-\dim I(d)=d+1\ge 2$
(see Lemma \ref{th:indetdim} for the computation $\dim I(d)=d$),
which yields $I(\Phi^{\mathrm{IRD}}_{d,n,k})=V_+(\Phi_n^*(V_{d^n,k}^{\vee}))=I(d)$ also by Lemma \ref{th:graded}.

Recall that $\overline{\rM}_{D,k}$, $D=d,d^n$,
is identified with $\Proj(\Sym(V_{D,k}^{\vee})^{\SL_{2,k}})$, and note
that the above $k$-algebra homomorphism $\Phi_n^*$ restricts to 
$\Phi_n^*:\Sym(V_{d^n,k}^{\vee})^{\SL_{2,k}}\to \Sym(V_{d,k}^{\vee})^{\SL_{2,k}}$. Recall the definition \eqref{eq:unstable}
 of $C_{d,n}$, and also \eqref{eq:ird} and \eqref{eq:coeff}.
Then we claim that
\begin{gather}
 V_+\bigl(\Phi_n^*
 \bigl((\Sym(V_{d^n,k}^{\vee})^{\SL_{2,k}})_+\bigr)\bigr)
 =q_{d,k}\Bigl(
 \bigl\{(\tilde{\phi}_{\xi_g},\dep_{\,\cdot}\tilde{\phi}_{\xi_g})
 \in I(d):\xi_g\in\MHRL(\phi)\bigr\}\cup C_{d,n}\Bigr);\label{eq:indetmod}
\end{gather}
indeed, for every $\phi\in\Rat_{d,K}$ and
every $F\in V_{d,k}$ such that $[F]=\hat{\phi}$,
we have $(\Phi_n^*s)(F)=s(F^{\circ n})$
for every $s\in(\Sym(V_{d^n,k}^{\vee})^{\SL_{2,k}})_+$, and
when $(\tilde{\phi}_{\xi_g},\dep_{\,\cdot}\tilde{\phi}_{\xi_g})\not\in I(d)$, recalling the definition (\cite[Definition 1.7 and Amplification 1.11]{MFK94}) of $\IRD_{d^n,k}^{\ss}(\xi_g)$ and \eqref{eq:git},
we have $(\Phi_n^*s)(F)=0$ for every $s\in(\Sym(V_{d^n,k}^{\vee})^{\SL_{2,k}})_+$ if and only if $\xi_g\notin\MHRL(\phi^n)$. Hence also by \eqref{eq:indet}, we have 
$q_{d,k}^{-1}(V_+(\Phi_n^*((\Sym(V_{d^n,k}^{\vee})^{\SL_{2,k}})_+)))
 =\{(\tilde{\phi}_{\xi_g},\dep_{\,\cdot}\tilde{\phi}_{\xi_g})
 \in I(d):\xi_g\in\MHRL(\phi)\}\cup C_{d,n}$.
When $d>2$, we have
$\dim q_{d,k}(\{(\tilde{\phi}_{\xi_g},\dep_{\,\cdot}\tilde{\phi}_{\xi_g})\in I(d):\xi_g\in\MHRL(\phi)\})\le d-3$,
and when $d=2$,
$\{(\tilde{\phi}_{\xi_g},\dep_{\,\cdot}\tilde{\phi}_{\xi_g})\in I(2):\xi_g\in\MHRL(\phi)\}=\emptyset$
(see Lemma \ref{th:indetdim}).
Once \eqref{eq:indetmod} is at our disposal, we have
$\operatorname{codim}_{\overline{\rM}_{d,k}}V_+(\Phi_n^*((\Sym(V_{d^n,k}^{\vee})^{\SL_{2,k}})_+))\ge(2d-2)-(d-2)\ge d\ge 2$
(see Theorem \ref{th:main}(\ref{item:moduli})
for the upper bound $d-2$ of $\dim q_{d,k}(C_{d,n})$),
which yields
$I(\Phi_{d,n,k})=V_+(\Phi_n^*((\Sym(V_{d^n,k}^{\vee})^{\SL_{2,k}})_+))=q_{d,k}(\{(\tilde{\phi}_{\xi_g},\dep_{\,\cdot}\tilde{\phi}_{\xi_g})\in I(d):\xi_g\in\MHRL(\phi)\}\cup C_{d,n})$ also by Lemma \ref{th:graded}.
\end{proof}

\begin{acknowledgement}
 This research was partially supported by JSPS Grant-in-Aid for Scientific Research (C), 23K03129, and
by the Research Institute for Mathematical Sciences,
an International Joint Usage/Research Center at Kyoto University.
\end{acknowledgement}


\begin{thebibliography}{99}

\bibitem{DeMarco}
Laura DeMarco.
\newblock Iteration at the boundary of the space of rational maps.
\newblock \emph{Duke Mathematical Journal} \textbf{130} (2005),
169--197.

\bibitem{DeMarcoJAMS}
Laura DeMarco.
\newblock The moduli space of quadratic rational maps.
\newblock \emph{Journal of the American Mathematical Society} \textbf{20} (2007),
321--355.

\bibitem{Faber}
Xander Faber.
\newblock Topology and geometry of the Berkovich ramification
locus for rational functions, I.
\newblock \emph{Manuscripta Mathematica} \textbf{142} (2013),
439--474.

\bibitem{Hartshorne}
Robin Hartshorne.
\newblock \emph{Algebraic Geometry}.
\newblock Graduate Texts in Mathematics, vol.~52, Springer-Verlag,
New York, 1977.

\bibitem{KN23}
Jan Kiwi and Hongming Nie.
\newblock Indeterminacy loci of iterate maps in moduli space.
\newblock {\em Indiana Univ. Math. J.}, 72(3):969--1026, 2023.

\bibitem{MFK94}
D.~Mumford, J.~Fogarty, and F.~Kirwan.
\newblock \emph{Geometric Invariant Theory}.
\newblock Third edition, Ergebnisse der Mathematik und ihrer Grenzgebiete (2),
vol.~34, Springer-Verlag, Berlin, 1994.

\bibitem{OkuyamaGeometry}
Y\^usuke Okuyama.
\newblock Geometric formulas on Rumely's weight function and crucial
measure in non-archimedean dynamics.
\newblock \emph{Mathematische Annalen} \textbf{376} (2020),
913--956.

\bibitem{OkuyamaDepth}
Y\^usuke Okuyama.
\newblock The intrinsic reductions and the intrinsic depths in
non-archimedean dynamics.
\newblock \emph{arXiv:2408.04415v3}.

\bibitem{OkuyamaQuadratic}
Y\^usuke Okuyama.
\newblock Intrinsic semistable reduction loci for the iterations of
non-archimedean quadratic rational functions and moduli-theoretic
applications.
\newblock \emph{arXiv:2603.03709v4}, 2026.

\bibitem{RumelyMinimal}
Robert Rumely.
\newblock The minimal resultant locus.
\newblock \emph{Acta Arithmetica} \textbf{169} (2015),
251--290.

\bibitem{RumelyEquivariant}
Robert Rumely.
\newblock A new equivariant in nonarchimedean dynamics.
\newblock \emph{Algebra \& Number Theory} \textbf{11} (2017),
841--884.

\bibitem{Seshadri77}
C.~S. Seshadri.
\newblock Geometric reductivity over arbitrary base.
\newblock \emph{Advances in Math.} \textbf{26} (1977),
225--274.

\bibitem{Silverman98}
Joseph~H. Silverman.
\newblock The space of rational maps on $\bP^1$.
\newblock \emph{Duke Math. J.} \textbf{94} (1998),
41--77.

\end{thebibliography}
\end{document}